\documentclass[a4paper,11pt]{article}
\usepackage[centertags]{amsmath}
\usepackage{amssymb,MnSymbol}
\usepackage{amsthm}
\usepackage{graphicx}
\usepackage[pagewise]{lineno}
 \usepackage{indentfirst}
\usepackage{amsfonts}
\usepackage{mathrsfs}
\usepackage{times,cite}
\usepackage{geometry}
\usepackage{color}
\newtheorem{thm}{Theorem}[section]
\newtheorem{lem}[thm]{Lemma}
\newtheorem{prop}[thm]{Proposition}

\theoremstyle{definition}
\newtheorem{rem}[thm]{Remark}

\newcommand{\e}{\varepsilon}
\newcommand{\p}{\partial}

\newcommand{\n}{\nonumber}
\newenvironment {Proof} {\noindent {\bf Proof. }}{\hfill $\blacksquare$\par\vspace{3mm}}

\numberwithin{equation}{section}
\allowdisplaybreaks    

\newenvironment{acknowledgment}{\noindent{\bf Acknowledgment }}{}
\date{}
\title{Rayleigh-Taylor instability for nonhomogeneous incompressible fluids with Navier-slip boundary conditions}
\author{{Shijin Ding$^a$, Zhijun Ji$^a$, Quanrong Li$^b$\thanks{Corresponding author.}}\\
{\it\small $^a$South China Research Center for Applied Mathematics and Interdisciplinary Studies,}\\
{\it\small South China Normal University,}{\it\small Guangzhou, 510631, Guangdong, China}\\
{\it\small $^b$College of Mathematics and Statistics,} {\it\small Shenzhen University,}\\
{\it\small Shenzhen, 518060, Guangdong, China}\\}

\begin{document}
\newcommand{\D}{\displaystyle}
\maketitle

 \begin{abstract}
 This paper is concerned with the {\it Rayleigh-Taylor} instability for the nonhomogeneous incompressible Navier-Stokes equations with Navier-slip boundary conditions around a steady-state in an infinite slab, where the Navier-slip coefficients do not have defined sign and the slab is horizontally periodic. Motivated by \cite{JiangJN}, we extend the result from Dirichlet boundary condition to Navier-slip boundary conditions. Our results indicate the factor that ``heavier density with increasing height'' still plays a key role in the instability under Navier-slip boundary conditions.
 \end{abstract}



\section{Introduction}
In this paper, we focus on the instability of the following nonhomogeneous incompressible Navier-Stokes equations with gravity in an infinite slab domain $\Omega=2\pi L\mathbb{T}\times(0,1)$:
\begin{equation}\label{1.1}
\begin{cases}
 \rho_t+\mathbf v\cdot \nabla\rho=0,\\
 \rho \mathbf v_t+\rho \mathbf v\cdot\nabla \mathbf v+\nabla p =\mu\Delta\mathbf v-\rho g \mathbf{e}_2 , \\
 {\rm div}\mathbf{v}=0.
\end{cases}
~~\text{in~~}\Omega,
\end{equation}
where the unknowns ($\rho,{\mathbf v},p$) denote density, velocity, and pressure of the fluid, respectively. The constant $\mu>0$ stands for the coefficient of shear viscosity, $\mathbf{e}_2=(0,1)$ is the vertical unit vector, and $-g\mathbf{e}_2$ describes the gravity.

The Navier-slip boundary conditions being considered is given as follows:
\begin{equation}\label{1.2}
\begin{cases}
\mathbf v\cdot \mathbf{n}=0,\\
2\mu \mathbb{D}(\mathbf{v})\cdot \mathbf{n}\cdot\mathbf{\tau}=k({\rm x})\mathbf v\cdot \tau ,
\end{cases}
~~\text{ on~~}\Sigma_1\bigcup\Sigma_0,
\end{equation}
where $\mathbb{D}({\mathbf v})=\frac{1}{2}(\nabla {\mathbf v}+\nabla^{T}{\mathbf v})$, $\mathbf{n}$ is the outward normal vector of the boundary and $\tau$ is the tangential vector, $\Sigma_1$ and $\Sigma_0$ are the upper and lower boundary, respectively, i.e. $\Sigma_1=2\pi L\mathbb{T}\times\{1\}, \Sigma_0=2\pi L\mathbb{T}\times\{0\}$, where $2\pi L\mathbb{T}$ stands for the 1D-torus of length $2\pi L$. In addition, $k({\rm x})$ is a scalar function describing the slip effect on the boundary. In this paper, $k({\rm x})$ will be taken to be constant $k_0$ and $k_1$ on $\Sigma_0$ and $\Sigma_1$, respectively, which do not have defined sign.

We first look for a smooth steady-state $(\bar\rho,0,\bar p)$ to the system (\ref{1.1}), where the density profile $\bar\rho:=\bar\rho(y)$ satisfies:
\begin{align}\label{1.3}
\bar\rho\in C^\infty([0,1]),\quad  \inf\limits_{y\in(0,1)}\bar\rho>0,
\end{align}
and the pressure $\bar p$ is determined by the following equality:
\begin{align}\label{1.4}
\nabla\bar p=-\bar\rho g\mathbf{e}_2.
\end{align}

Since we are interested in Rayleigh-Taylor(RT) instability, we assume that the steady density satisfies
\begin{align}\label{1.5}
   {\bar\rho'(y_0)}>0,\quad \text{for~some~} y_0\in(0,1).
\end{align}
This condition means that there is a neighborhood of $y_0$, such that $\bar{\rho}$ increases with $y$, i.e., a heavy
fluid is on top of the light one.

Now, we define the perturbation as
\begin{align*}
(\varrho,\mathbf u,q):=(\rho-\overline\rho,\mathbf v-\mathbf 0,p-\overline p).
\end{align*}
Then, $(\varrho,\mathbf u,q)$ satisfies the following equations
\begin{equation}\label{1.6}
\begin{cases}
\varrho_t+\mathbf u\cdot\nabla(\varrho+\bar\rho)=0,\\
(\varrho+\bar\rho)\mathbf u_t+(\varrho+\bar\rho){\mathbf u}\cdot\nabla{\mathbf u}+\nabla q+\varrho g\mathbf{e}_2=\mu\Delta {\mathbf u},\\
\mathrm{div}{\mathbf u}=0,
\end{cases}
~\text{in~}\Omega,
\end{equation}
with the corresponding initial data and boundary conditions turning to
\begin{align}\label{1.7}
\begin{cases}
(\varrho,\mathbf u)|_{t=0}=(\varrho_0,\mathbf u_0),&\text{in~~} \Omega,\\
u_2=0,\ &\text{on~~}\Sigma_1\bigcup\Sigma_0,\\
\p_yu_1=\frac{k_1}{\mu}u_1,\ &\text{on~~}\Sigma_1,\\
\p_yu_1=-\frac{k_0}{\mu}u_1,\ &\text{on~~}\Sigma_0,.
\end{cases}
\end{align}

Linearizing system (\ref{1.6}), one gets
\begin{equation}\label{1.8}
\begin{cases}
\varrho_t+\bar\rho'u_2=0,\\
\bar\rho\mathbf u_t+\nabla q+\varrho g\mathbf{e}_2=\mu\Delta\mathbf u,\\
{\rm div}{\mathbf u}=0.\\
\end{cases}
\text{in~~}\Omega.
\end{equation}

To analyze our problem, we would like to apply the growing normal mode method, for which the readers can refer to \cite{S.Ding}, for instance.

Precisely, we first assume a growing mode ansatz of solutions
\begin{align*}
(\varrho,{\mathbf u},q)(x,y;t)=e^{\lambda t}({\tilde\rho},\tilde{\mathbf v},{\tilde p})(x,y)
\end{align*}
to the linearized system (\ref{1.8}) with some constant $\lambda >0$. Substituting these ansatz into \eqref{1.8} deduces a system for the new unknowns $({\tilde\rho},\tilde{\mathbf v},{\tilde p})$ as follows
\begin{equation}\label{1.9}
\begin{cases}
\lambda\tilde\rho+\bar\rho'\tilde v_2=0,\\
\lambda\bar\rho\tilde{\mathbf v}+\nabla{\tilde p}+\tilde\rho g\mathbf{e}_2=\mu\Delta \tilde{\mathbf v},\\
{\rm div}\tilde {\mathbf v}=0.\\
\end{cases}
~\text{in~}\Omega.
\end{equation}
Eliminating $\tilde\rho$ in \eqref{1.9}, one arrives at
\begin{align}\label{1.10}
\begin{cases}
{\lambda}^2\bar\rho\tilde{\mathbf v}+\lambda\nabla{\tilde p}=\lambda\mu\Delta \tilde{\mathbf v}+g\bar\rho' {\tilde v}_2\mathbf{e}_2,\\
\text{div}\tilde{\mathbf v}=0,\\
\end{cases}
\end{align}
which is endowed with the following initial data and boundary conditions:
\begin{align}\label{1.11}
\begin{cases}
(\tilde\rho,\tilde {\mathbf v})|_{t=0}=(\varrho_0,\mathbf u_0),&\text{in~~} \Omega,\\
\tilde v_2=0,\ &\text{on~~}\Sigma_1\bigcup\Sigma_0,\\
\p_y{\tilde v}_1=\frac{k_1}{\mu}{\tilde v}_1,\ &\text{on~~}\Sigma_1,\\
\p_y{\tilde v}_1=-\frac{k_0}{\mu}{\tilde v}_1,\ &\text{on~~}\Sigma_0.
\end{cases}
\end{align}

Second, for any frequency $\xi\in\mathbb{R}$, $\xi\neq 0$, we rewrite the unknowns in \eqref{1.10}-\eqref{1.11} in terms of $(\phi,\psi,\pi)(y):(0,1)\rightarrow\mathbb{R}$ with
\begin{align*}
\begin{cases}
\widetilde v_1(x,y)=-i\phi(y)e^{ix\xi},\\
\widetilde v_2(x,y)=\psi(y)e^{ix\xi},\\
\widetilde p(x,y)=\pi(y)e^{ix\xi},
\end{cases}
\end{align*}
from which we can infer that the new unknowns $(\phi,\psi,\pi)$ satisfy the following ODEs
\begin{align}\label{1.12}
\begin{cases}
-{\lambda}^2\bar\rho\phi+\lambda\xi\pi-\lambda\mu({\xi}^2\phi-\phi'')=0,\\
{\lambda}^2\bar\rho\psi+\lambda\pi'+\lambda\mu({\xi}^2\psi-\psi'')=g\bar\rho'\psi,\\
\xi\phi+\psi'=0,\\
\end{cases}
\end{align}
with the corresponding boundary conditions
\begin{align}\label{1.13}
\begin{cases}
\psi(0)=\psi(1)=0,\\
\phi'(1)=\frac{k_1}{\mu}\phi(1),\\
\phi'(0)=-\frac{k_0}{\mu}\phi(0).\\
\end{cases}
\end{align}

Third, eliminating $\pi$ in \eqref{1.12}, one yields a fourth-order ODE of $\psi(y)$:
\begin{align}\label{1.14}
-\lambda^2 \left[\xi^2\bar\rho\psi-\left(\bar\rho\psi'\right)' \right]
=\lambda\mu\left(\psi^{(4)}-2\xi^2\psi''+\xi^4\psi\right)-g\xi^2\bar\rho'\psi,
\end{align}
and the corresponding boundary conditions：
\begin{align}\label{1.15}
\begin{cases}
\psi(0)=\psi(1)=0,\\
\psi''(1)=\frac{k_1}{\mu}\psi'(1),\\
\psi''(0)=-\frac{k_0}{\mu}\psi'(0).\\
\end{cases}
\end{align}

In conclusion, the problem \eqref{1.7}-\eqref{1.8} is finally reduced to the fourth-order ODE system \eqref{1.14}-\eqref{1.15}.

The main results of this paper are stated as follows:
\begin{thm}\label{th1}
(Linear instability) The steady-state $(\bar\rho,0,\bar p)$ is linearly unstable, provied that the steady density profile $\bar\rho$ satisfies \eqref{1.3} and \eqref{1.5}. Precisely, there exists exponentially growing solution to the linearized perturbed problem \eqref{1.7}-\eqref{1.8}, such that $\rVert (\varrho,\mathbf{u},q)(t)\rVert_{H^k(\Omega)}\rightarrow\infty$, as $t\rightarrow\infty$.
\end{thm}

\begin{thm}\label{th2}
(Nonlinear instability) The steady state solution $(\bar\rho,0,\bar p)$ is nonlinearly unstable in the Hardama sense, provided that the steady density profile $\bar\rho$ satisfies \eqref{1.3} and \eqref{1.5}. Precisely, there exist positive constants $\Lambda$, $\e$, $m_0$, and a pair $(\bar{\varrho}_0,\bar{\mathbf{u}}_0)\in L^2(\Omega)\times H^2(\Omega)$, such that for any $\delta\in (0,\e)$, there is an unique global strong solution $({\varrho},\mathbf{u})$ to the nonlinear perturbed problem \eqref{1.6}-\eqref{1.7} on $[0,T)$ with the initial data
 $(\varrho_0,\mathbf{u}_0):=(\delta\bar{\varrho}_0,\delta\bar{\mathbf{u}}_0)$, but
\begin{equation}\label{1.16}
\|\varrho(T^\delta)\|_{L^2(\Omega)},\ \|\mathbf{u}(T^\delta)\|_{L^2(\Omega)}\geq {\e},
\end{equation}
for some escape time $T^\delta:=\frac{1}{\Lambda}\ln\frac{2\e}{m_0\delta}\in(0,T)$.
\end{thm}

Before proving these theorems, let us recall some results on the RT instability problems. The RT instability is a kind of well-known instability in fluid dynamics, which is driven by the gravity when the upper fluid is heavier than the lower one. In 1883, Rayleigh {\cite{Rayleigh}} first considered the linear instability for an incompressible fluid. 120 years later, in 2003, Hwang and Guo {\cite{Hwang.Guo}} studied 2D nonhomogeneous incompressible inviscid fluid in strip domain with zero normal velocity on the boundary, and proved the RT instability in the Hardama sense. However, when the viscosity is taken into account, there is no direct variational structure for constructing exponentially growing solution. In 2011, Guo and Tice {\cite{Guo.Tice11}} introduced a general method for these problem, they turned to study the modified variational problem first and then went back to the original problem by fixed-point theory. Motivated by {\cite{Guo.Tice11}}, F. Jiang, S. Jiang and G. Ni {\cite{JiangJN}} in 2013 studied the RT instability of nonhomogeneous incompressible viscous fluid at the present of the uniform gravity field in ${\mathbb{R}^3}$. Whereafter, F. Jiang and S. Jiang made a breakthrough in three-dimensional bounded domain case{\cite{JiangJ}} and analyzed both the instability and stability for given different steady density profile $\bar\rho$.

When fluids are electrically conducting at the present of magnetic field, the RT instability arises and the growth of the instability will be influenced by the magnetic field due to the Lorentz force. Thus, there are also some authors pay considerable attention to the inhibition of the magnetic field on the RT instability. In 1954,  Kruskal and Schwarzchild first proved that a horizontal magnetic field has no effect on the development of the linear RT instability \cite{Kruskal}. Afterwards, Hide in \cite{Hide} investigated the influence of a vertical magnetic field. To our knowledge, there is some critical magnetic number $B_c$, such that when the vertical background magnetic field is less than the critical magnetic number, then the magnetic field has no effect on the RT instability \cite{JJW}, while the vertical background magnetic field is larger than the critical magnetic number, then the magnetic field has an inhibitory on the RT instability \cite{WYJ}. Therefore, this physical phenomenon has been verified mathematically in some cases.

In 2018, Ding, Li and Xin {\cite{S.Ding}} found that there exists a critical viscosity coefficient for distinguishing stability from instability when considering the homogeneous fluid with Navier-slip boundary conditions in 2D strip domain.

As for the Navier-slip boundary conditions, which was first proposed by C. Navier {\cite{C.Navier}} in 1827, it describes the phenomenon that fluid moves along the boundary. Mathematically, it can be depicted as in \eqref{1.2}. Compared to the classical Dirichlet boundary condition, Navier-slip boundary conditions are more realistic in some situations, but the mathematical literature is less. Early in 1973, V. Solonnikov and V. \v{S}\v{c}ailov \cite{V.S} gave the first rigorous mathematical analysis to Navier-Stokes equations with Navier-slip boundary conditions, they focused on linear stationary equations in 3D with the coefficient $k({\rm x})=0$ and the external force $\mathbf{f}\in L^2$. In 1980s, G. Mulone and F. Salemi \cite{M.F83,M.F85} considered the Navier-Stokes equations with Navier-slip boundary conditions in a three dimensional bounded domain. They proved the well-posedness for the corresponding stationary problem, and the existence of weak solution for the evolutionary problem. Afterwards, J.Kelliher {\cite{Kelliher}} established the existence, uniqueness and regularity for 2D bounded domain case in 2006 , where the domain is consisted of a finite number of connected components with the slip coefficient $k({\rm x}) \in L^\infty(\p\Omega)$. Based on the above results of weak solution, in 2010, H. B. da Veiga {\cite{HB.Veiga}} improved the regularity of the weak solutions, up to the boundary. Later, C. Amrouche and co-authors studied stationary and evolutionary problems{\cite{C.A14,C.A16,C.A11}} in $L^p$ with $p\in(1,\infty)$, proving the existence of weak and strong solutions in a three-dimensional bounded domain with smooth boundary. It should be noted that, the noncompact infinite slab case is not included in any result mentioned above. In 2018, Ding and Li proved the existence and uniqueness of strong solution of incompressible fluid with Navier boundary conditions in 2D infinite slab \cite{Li}. For the references on the vanishing viscosity limit of Navier-Stokes equation with Navier boundary conditions, we refer the readers to {\cite{T.Clo,XIAO.XIN07,XIAO.XIN13}} and the reference therein. The main aim of this paper is to investigate the slip effect from the boundary on the RT instability.

The rest of the paper is arranged as follows. In Section 2, we give some notations and list some useful inequalities. In Section 3 and Section 5, we will give the proofs of Theorem \ref{th1} and Theorem \ref{th2}, respectively. In Section 4, we deduce the energy estimates, which is a preparation for proving the nonlinear instability in Section 5.

\section{Preliminary}

For simplicity, we denote $L^2(0,1)$ and $H^k(0,1)$ by $L^2$ and $H^k$. Without confusion, we will also write $L^p(\Omega)$ and $H^k(\Omega)$  by $L^p$ and $H^k$, respectively. The integral form $ \int_\Omega fdxdy$ will be simply denoted by $\int f$. In addition, the scalar function and vector function will be denoted by $f$ and $\mathbf{f}$ for distinction, such as $\mathbf{f}=(f_1,f_2) $. The product functional space $(X)^2$ will also be denoted by $X $, for example, the vector function $\mathbf{u}\in (H^1)^2$ will be still denoted by $\mathbf{u}\in H^1$. The usual notations will be used as in general unless with extra statements.

For convenience, we list a few lemmas that will be used in this paper, without the proofs. The readers interested in the proof could refer to \cite{Li} for the details.
\begin{lem}\label{lemma2.1}
{\rm (Poincar\'{e} inequality in $\Omega$)} There exists constant $C>0$, such that for any $\mathbf{u}\in V=\{\mathbf{v}\in H^1|\nabla\cdot{\mathbf{v}}=0, {\mathbf{v}}\cdot{\mathbf{n}}=0\textrm {~on~} \p\Omega \}$, the following inequality holds:
\begin{equation}\label{2.1}
   \|\mathbf{u}\|_{L^2(\Omega)}\leq C\|\p_y\mathbf{u}\|_{L^2(\Omega)}.
\end{equation}
\end{lem}

\begin{lem}\label{lemma2.2}
{\rm ($L^4$ estimate in $\Omega$)} There exists constant $C>0$, such that for any $\mathbf{u}\in W=\{\mathbf{v}\in V\cap H^2|\mathbf{v} {\textrm ~satisfies ~} (\ref{1.2})_2\},$
there holds
\begin{equation}\label{2.2}
\|\mathbf{u}\|_{L^4(\Omega)}^2\leq C\|\mathbf{u}\|_{L^2(\Omega)}\|\nabla\mathbf{u}\|_{L^2(\Omega)}.
\end{equation}
\end{lem}

\begin{lem}\label{lemma2.3}
{\rm($L^\infty$ estimate in $\Omega$)} There exists constant $C>0$, such that for any $\mathbf{u}\in W$ , one has
\begin{equation}\label{2.3}
\|\mathbf{u}\|_{L^\infty(\Omega)}^2\leq C\|\mathbf{u}\|_{L^2(\Omega)}\|\mathbf{u}\|_{H^2(\Omega)}.
\end{equation}
\end{lem}

\begin{lem}\label{lemma2.4}
{\rm($L^2$ estimate of gradient in $\Omega$)} There exists constant $C>0$, such that for any $\mathbf{u}\in W$, there holds
\begin{equation}\label{2.4}
\|\nabla\mathbf{u}\|_{L^2(\Omega)}^2\leq C\|\mathbf{u}\|_{L^2(\Omega)}\|\nabla\mathbf{u}\|_{H^1(\Omega)}.
\end{equation}
\end{lem}

\begin{lem}\label{lemma2.5}
{\rm($L^4$ estimate of gradient in $\Omega$)} There exists constant $C>0$, such that for any $\mathbf{u}\in W$, there holds
\begin{equation}\label{2.5}
\|\nabla\mathbf{u}\|_{L^4(\Omega)}^2\leq C\|\mathbf{u}\|^{\frac{1}{2}}_{L^2(\Omega)}\|\nabla\mathbf{u}\|^{\frac{3}{2}}_{H^1(\Omega)}.
\end{equation}
\end{lem}
\Proof Using estimate \eqref{2.3}, the proof of \eqref{2.5} is similar to that of \eqref{2.4}.\endProof
\section{The proof of Theorem \ref{th1}}
In this section, we will make effort to construct a solution for (\ref{1.14})-(\ref{1.15}) by variational method, which at once deduces a solution for ODE system (\ref{1.12})-(\ref{1.13}). Then, an exponentially growing solution will be given for linearized perturbed equation (\ref{1.8})-(\ref{1.7}), and Theorem\ref{th1} follows.

However, similar to \cite{Guo.Tice11}, the appearance of $\lambda$ both quadratically and linearly in equation \eqref{1.14} breaks the natural variational structure of problem \eqref{1.14}-\eqref{1.15}. To circumvent this obstacle, for any spatial frequency $\xi$ ($\neq 0$), we introduce a family $(s>0)$ of modified variational problem
\begin{align}\label{3.1}
 \alpha(s):=\inf\limits_{\psi\in\mathcal{A}}E(\psi,s),
\end{align}
with the energy functional
\begin{align}\label{3.2}
E(\psi,s):=&s \left[\int_0^1\mu|\psi''|^2dy-(k_1|\psi'(1)|^2+k_0|\psi'(0)|^2)\right]\nonumber\\
 &+s\mu\int_0^1\left(2\xi^2|\psi'|^2+\xi^4\psi^2\right)dy-\int_0^1 g\xi^2\bar\rho'\psi^2dy,
\end{align}
and the corresponding admissible set
\begin{align}\label{3.3}
\mathcal{A}=\left\{\psi\in H_0^1\cap H^2\Big|J(\psi):=\int_0^1\bar\rho\left(\xi^2\psi^2+|\psi'|^2 \right)dy=1 \right\}.
\end{align}
In addition, to emphasize the dependence on $\xi$, we sometimes write \eqref{3.1} as
\begin{align*}
\alpha(s,\xi)=\inf\limits_{\psi\in\mathcal{A}}E(\psi,s,\xi).
\end{align*}

In order to recover the corresponding variational form of problem \eqref{1.14}-\eqref{1.15}, we will show later in this section that there exists a fixed point $s_0>0$, so that $\alpha(s_0)=-s_0^2$. The first proposition is devoted to proving the well-definedness of \eqref{3.1}.
\begin{prop}\label{pr3.1}
Suppose that the steady density satisfies (\ref{1.3}). Then, for any $\xi\neq 0$ and $s\in(0,+\infty)$,
$E(\psi,s)$ achieves its minimum on $\mathcal{A}$.
\end{prop}
\Proof
For any $\psi\in\mathcal{A}$, since $J(\psi)=1$, by Cauchy inequality, one has

\begin{align}\label{3.4}
E(\psi)
\geq& -\int_0^1 g\xi^2\bar\rho'\psi^2dy+s \int_0^1\{\mu|\psi''|^2-[((k_0+k_1)y-k_0)|\psi'(y)|^2]_y\}dy\nonumber\\
\geq&-g\left\|\frac{\bar\rho'}{\bar\rho}\right\|_{L^{\infty}}+s \int_0^1\{\mu|\psi''|^2-\sum_{i=0}^1k_i|\psi'|^2-2\sum_{i=0}^1k_iy-k_0)\psi'\psi''\}dy\nonumber\\
\geq& -g\left\|\frac{\bar\rho'}{\bar\rho}\right\|_{L^{\infty}}-sC_0\int_0^1|\psi'|^2dy
\geq -g\left\|\frac{\bar\rho'}{\bar\rho}\right\|_{L^{\infty}}-sC_0\left\|\frac{1}{\bar\rho}\right\|_{L^{\infty}},
\end{align}

where $C_0:=\max\limits_{0\leq y\leq1}[|k_0+k_1|+\mu^{-1}((k_0+k_1)y-k_0)^2]$.

Since $\bar\rho$ has a positive lower bound, it follows from \eqref{3.4} that $E(\psi)$ is bounded below on $\mathcal{A}$, and so it has an infimum. Denote that
\begin{align*}
M:=\inf\limits_{\psi\in\mathcal{A}}E(\psi).
\end{align*}
In view of the definition of infimum, there exists a  minimizing sequence $\{\psi_n\}_{n=1}^\infty\in\mathcal{A}$, such that
\begin{align*}
{\lim_{n\to+\infty}}E(\psi_n)=M.
\end{align*}
Without loss of generality, we suppose that
\begin{align}\label{3.5}
  E(\psi_n)\leq M+1.
\end{align}

On one hand, since $J(\psi_n)=1$, there holds
\begin{align}\label{3.6}
\int_0^1\psi_n^2dy\leq\xi^{-2}\left\|\frac{1}{\bar\rho}\right\|_{L^\infty} \quad\text{and}\quad\int_0^1|\psi'_n|^2dy\leq\left\|\frac{1}{\bar\rho}\right\|_{L^\infty},
\end{align}
which imply that $\{\psi_n\}_{n=1}^\infty$ is a bounded sequence in $H_0^1$.

On the other hand, similar to (\ref{3.4}), it follows from the definition of $E(\psi)$ that
\begin{align}\label{3.7}
  \mu\int_0^1|\psi_n''|^2dy
  &=\frac{E(\psi_n)+\int_0^1g\xi^2\bar\rho'\psi_n^2dy}{s}+\sum_{i=0}^1k_i|\psi_n'(i)|^2-\mu\xi^2\int_0^1(2|\psi_n'|^2+\xi^2\psi_n^2)dy\nonumber\\
  &\leq\frac{M+1+g\left\|\frac{\bar\rho'}{\bar\rho}\right\|_{L^\infty}}{s}+2C_0\int_0^1|\psi_n'|^2dy+\frac{\mu}{2}\int_0^1|\psi_n''|^2dy,
\end{align}
which, together with (\ref{3.6}), indicate that
\begin{align}\label{3.8}
\int_0^1|\psi_n''|^2dy
\leq2\mu^{-1}\left[s^{-1}\left(M+1+g\left\|\frac{\bar\rho'}{\bar\rho}\right\|_{L^\infty}\right)+2C_0\left\|\frac{1}{\bar\rho}\right\|_{L^\infty}\right].
\end{align}
 In conclusion, for any fixed $s\in(0,\infty)$, $\{\psi_n\}_{n=1}^\infty$ is a bounded sequence in $H_0^1\cap H^2$, and hence, there exists a convergent subsequence, still denoted by $\{\psi_{n}\}_{n=1}^\infty$, such that $\psi_{n}\rightharpoonup\psi$ weakly in $H^2$ and $\psi_{n}\rightarrow\psi$ strongly in $H_0^1$. Then $J(\psi)=\lim\limits_{n\to+\infty}J(\psi_{n})=1$, which implies $\psi\in\mathcal{A}$.

In addition, using the weak lower semi-continuity and the convergence of $\{\psi_{n}\}_{n=1}^\infty$, one gets
\begin{align}\label{3.9}
  E(\psi)
  =~&s\mu\int_0^1|\psi''|^2dy+s(2\mu\xi^2-(k_1+k_0))\int_0^1|\psi'|^2dy\nonumber\\
   ~&+\xi^2\int_0^1{(\xi^2s\mu-g\bar\rho')\psi^2}dy-s\int_0^1((k_1+k_0)y-k_0)\psi'\psi''dy\nonumber\\
  \leq~&{\liminf\limits_{n\to+\infty}}E(\psi_{n})=\lim\limits_{n'\rightarrow\infty}E(\psi_{n'})=M,
\end{align}
 which means that $\psi$ is the minimizer of $E(\cdot)$ on $\mathcal{A}$. The proof of this proposition is finished.
\endProof

In what follows, we will prove that $\alpha(s)$ is negative for any $s \in (0,\mathfrak{S})$, where $\mathfrak{S}$ will be defined in Proposition \ref{pr3.4}. For convenience, we denote
\begin{eqnarray*}
\begin{cases}
  E_0(\psi)=\int_0^1\mu|\psi''|^2dy-(k_1|\psi'(1)|^2+k_0|\psi'(0)|^2), \\
  E_1(\psi)=\mu\int_0^1(2|\psi'|^2+\xi^2\psi^2)dy, \\
  G(\psi)=E_0(\psi)+\xi^2E_1(\psi),\\
  E_2(\psi)=\int_0^1g\xi^2\bar\rho'\psi^2dy.
\end{cases}
\end{eqnarray*}
Then, \eqref{3.2} can be rewritten as $E(\psi)=sG(\psi)-E_2(\psi)$.

Similar to \cite{S.Ding}, firstly, we define the critical viscosity:
\[\mu_c:=\sup\limits_{\psi\in\mathcal{A}}\frac{k_1|\psi'(1)|^2+k_0|\psi'(0)|^2}{\int_0^1|\psi''|^2dy},\]
which, according to \cite{S.Ding}, can be explicitly expressed by
\begin{equation*}
\mu_c=
\begin{cases}
0,&~k_0\leq 0~\&~k_1\leq 0;\\
k/6,&~k_0=k_1:=k>0;\\
\frac{(k_0+k_1)+\sqrt{k_0^2+k_1^2-k_0k_1}}{6},&~\text{otherwise}.
\end{cases}
\end{equation*}
It is clear that $E_0(\psi)\geq0$ for any $\psi\in\mathcal{A}$, provided $\mu\geq\mu_c$. Conversely, there exists $\psi\in\mathcal{A}$ such that $E_0(\psi)<0$, if $\mu<\mu_c$.

Further, according to \cite{S.Ding}, if $\mu<\mu_c$, there exists a critical frequency defined by
\[0<\xi_c^2:=\sup\limits_{\psi\in\mathcal{A}}\frac{k_1|\psi'(1)|^2+k_0|\psi'(0)|^2-\int_0^1\mu|\psi''|^2dy}{\mu\int_0^1(2|\psi'|^2+\xi_c^2\psi^2)dy}
  =\sup\limits_{\psi\in\mathcal{A}}\frac{-E_0(\psi)}{E_1(\psi,\xi_c)}.\]
And, under the assumption of $\mu<\mu_c$, one can point out from the definition of $\xi_c$ that
\begin{align*}
  0<\xi_c^2\leq\frac{k_1|\psi'(1)|^2+k_0|\psi'(0)|^2-\int_0^1\mu|\psi''|^2dy}{\mu\int_0^1 2|\psi'|^2dy}
  \leq\frac{C_0\int_0^1|\psi'|^2dy}{2\mu\int_0^1|\psi'|^2dy}=\frac{C_0}{2\mu},
\end{align*}
which indicates that $0<|\xi_c|\leq\sqrt{\frac{C_0}{2\mu}}$. Here, $C_0$ is the same number defined as in (\ref{3.4}). In conclusion, for any $\mu>0$, the functional $G(\psi)$ is positive on $\mathcal{A}$, provided $|\xi|>\xi_c$.

Now, we define the solvable domain for frequency
 \[\mathbb{A}^g:=\{\xi\in \mathbb{R}||\xi|\in(a,b):=(|\xi_c|,b)\subset(0,+\infty)\},\]
 where $0<a<<b<+\infty$. Next, we will show that for any $\xi\in\mathbb{A}^g$, there exists a fixed point $s_0>0$ such that $\alpha(s_0)=-s^2_0$. To this end, several propositions will be devoted to study some properties of $\alpha(s)$ in what follows.

\begin{prop}\label{pr3.2}
For any $\xi\in\mathbb{A}^g$ and $s\in(0,+\infty)$, there exist positive numbers $C_1,C_2$, depending on $\mu,a,b,\bar\rho,g,k_0,k_1$, such that $\alpha(s)\leq sC_2-C_1$.
\end{prop}
\Proof
In view of \eqref{1.5}, there is a point $y_0\in (0,1)$, such that $\bar\rho'(y_0)>0$. Thus, there exists an open neighbourhood of $y_0$
\begin{align*}
B_{y_0}^{\delta}:=\left\{y\in(0,1)\big\vert~|y-y_0|<\delta\right\}\subset(0,1)
\end{align*}
with $\delta$ small sufficiently, such that $\bar\rho'(y)>\bar\rho'(y_0)/2$ in $B_{y_0}^{\delta}$. Now, define
\[\tilde\psi(y):=f(y-y_0),\]
where $f(r)\in C_0^{\infty}(\mathbb{R})$ is a cut-off function satisfying $f(r)=1$ for $|r|<\delta/2$ and $f(r)=0$ for $|r|\geq3\delta/4$. Then, $\tilde\psi(y)\in H_0^1\cap H^2$ and $E_2(\tilde\psi)=\int_0^1{g\xi^2\bar\rho'\tilde\psi^2}dy>g\xi^2\bar\rho'(y_0)\delta/4>0$.

So, for any $\xi\in\mathbb{A}^g$, there holds that
\begin{align}\label{3.10}
  \frac{E_2(\tilde\psi)}{J(\tilde\psi)}
  =\frac{g\xi^2\int_0^1\bar\rho'\tilde\psi^2dy}{\int_0^1\bar\rho(\xi^2\tilde\psi^2+|\tilde\psi'|^2)dy}
  \geq\frac{ga^2\int_0^1\bar\rho'\tilde\psi^2dy}{\int_0^1\bar\rho(b^2\tilde\psi^2+|\tilde\psi'|^2)dy}:=C_1,
\end{align}
and that
\begin{align}\label{3.11}
  \frac{G(\tilde\psi)}{J(\tilde\psi)}
   & = \frac{\mu\int_0^1(|\tilde\psi''|^2+2\xi^2|\tilde\psi'|^2+\xi^4\tilde\psi^2)dy-(k_1|\tilde\psi'(1)|^2+k_0|\tilde\psi'(0)|^2)}
   {\int_0^1\bar\rho(\xi^2\tilde\psi^2+|\tilde\psi'|^2)dy}\nonumber\\
   &\leq\frac{\mu\int_0^1(2b^2|\tilde\psi'|^2+b^4\tilde\psi^2)dy+C_0\left\|\frac{1}{\bar\rho}\right\|_{L^{\infty}(0,1)}}
      {\int_0^1\bar\rho(a^2\tilde\psi^2+|\tilde\psi'|^2)dy}
   := C_2,
\end{align}
where $C_1$ and $C_2$ are positive constants depending on $a,b,\bar\rho,\mu,k_0,k_1$.

In conclusion, recalling that
\begin{align}\label{3.12}
\alpha(s)
&=\inf\limits_{\substack{\psi\in H_0^1\cap H^2\\ \psi\neq 0}} \frac{E(\psi)}{J(\psi)}
:=\frac{sG(\psi)-E_2(\psi)}{\int_0^1\bar\rho(\xi^2\psi^2+|\psi'|^2)dy},
\end{align}
one yields
\begin{align}\label{3.13}
  \alpha(s)
  =\inf\limits_{\substack{\psi\in H_0^1\cap H^2\\ \psi\neq 0}} \frac{E(\psi)}{J(\psi)}
  \leq\frac{E(\tilde\psi)}{J(\tilde\psi)}\leq sC_2-C_1,
\end{align}
for any $\xi\in\mathbb{A}^g$.
\endProof

\begin{prop}\label{pr3.3}
For any $\xi\in\mathbb{A}^g$, $\alpha(s)$ is a locally Lipschitz continuous function on $(0,\infty)$, that is, $\alpha(s)\in C_{\rm{loc}}^{0,1}(0,\infty)$.
\end{prop}
\Proof
For any $[c,d]\subset(0,\infty)$ and $s\in[c,d]$, in view of Proposition\ref{pr3.1}, there exists a minimizer $\psi_s\in\mathcal{A}$, such that $\alpha(s)=E(\psi_s,s)$. Then, by the definition of $G(\psi_s,s)$, one has
\begin{align}\label{3.14}
  G(\psi_s,s)
  = & \frac{E(\psi_s,s)}{s}+\frac{g}{s}\int_0^1\xi^2\bar\rho'|\psi_s|^2dy
  \leq \frac{E(\psi_s,s)}{s}+\frac{g}{c}\left\|\frac{\bar\rho'}{\bar\rho}\right\|_{L^{\infty}}\n\\
  \leq & \frac{1+dC_2-C_1}{c}
  +\frac{g}{c}\left\|\frac{\bar\rho'}{\bar\rho}\right\|_{L^{\infty}} :=K,
\end{align}
where $K$ is a positive number depending on $c,d,\mu,\bar\rho,g,k_0,k_1$.

Similarly, for any $s_1,s_2\in[c,d]$, there exist $\psi_{s_1}, \psi_{s_2}\in\mathcal{A}$, such that $\alpha(s_1)=E(\psi_{s_1},s_1)$ and $\alpha(s_2)=E(\psi_{s_2},s_2)$. Furthermore, note that
\begin{align*}
  \alpha(s_1)&=E(\psi_{s_1},s_1)\leq E(\psi_{s_2},s_1)=E(\psi_{s_2},s_2)+(s_1-s_2)G(\psi_{s_2})\leq\alpha(s_2)+K|s_1-s_2|,\\
\alpha(s_2)&=E(\psi_{s_2},s_2)\leq E(\psi_{s_1},s_2)=E(\psi_{s_1},s_1)+(s_2-s_1)G(\psi_{s_1})\leq\alpha(s_1)+K|s_1-s_2|.
\end{align*}
Hence, one deduces
\begin{equation}\label{3.15}
  |\alpha(s_1)-\alpha(s_2)|\leq K|s_1-s_2|,
\end{equation}
which implies that $\alpha(s)\in C_{\rm{loc}}^{0,1}(0,\infty)$.
\endProof

By virtue of Proposition \ref{pr3.2}, there exists a constant $s^*>0$ such that
\begin{equation}\label{3.16}
  \alpha(s)<0~, \text{~for~any~} s\in(0,s^*),
\end{equation}
which implies the following proposition

\begin{prop}\label{pr3.4}
For any $\xi\in\mathbb{A}^g$, denote $\mathfrak{S} :=\inf\{\tau~|~\alpha(\tau)>0\}$. Then $\alpha(s)<0$ holds for any $s\in(0,\mathfrak{S})$.
\end{prop}


Next, we will prove the well-definedness of $-\lambda^2=\alpha(\lambda)$ with $\lambda\in(0,\mathfrak{S})$ by intermediary theorem.
\begin{prop}\label{pr3.5}
For any $\xi\in\mathbb{A}^g$, there exists $\lambda \in(0,\mathfrak{S})$, such that
\begin{align}\label{3.17}
-\lambda^2(\xi)=\alpha(\lambda(\xi))=\inf\limits_{\psi\in\mathcal{A}}E(\psi,\lambda(\xi)).
\end{align}
\end{prop}
\Proof
For any $s_1,s_2\in(0,\mathfrak{S})$, $s_1<s_2$, and denote that $\psi_{s_i}$ is the minimizer of functional $E(\psi,s_i)$ on $\mathcal{A},~i=1,2$, Then, it follows from $G(\cdot)>0$ that
\begin{equation}\label{3.18}
  \alpha(s_1)=E(\psi_{s_1},s_1)\leq E(\psi_{s_2},s_1)<E(\psi_{s_2},s_2)=\alpha(s_2),
\end{equation}
This indicates that $\alpha(s)$ is strictly monotonically increasing on $s\in(0,\mathfrak{S})$.

Now, define
\begin{equation}\label{3.19}
  \Phi(s)=\frac{s^2}{-\alpha(s)}.
\end{equation}

On one hand, by proposition \ref{pr3.2}, one has
\begin{equation}\label{3.20}
  \lim\limits_{s\rightarrow 0^+}-\alpha(s)\geq\lim\limits_{s\rightarrow 0^+}-(sC_2-C_1)=C_1>0,
\end{equation}
which implies that
\begin{equation}\label{3.21}
  \lim\limits_{s\rightarrow 0^+}\Phi(s)
  =\lim\limits_{s\rightarrow 0^+}\frac{s^2}{-\alpha(s)}
  =0.
\end{equation}

On the other hand, from the definition of $\mathfrak{S}$, it is clear that
\begin{equation}\label{3.22}
  \lim\limits_{s\rightarrow\mathfrak{S}^-}\Phi(s)=\lim\limits_{s\rightarrow\mathfrak{S}^-}\frac{s^2}{-\alpha(s)}=+\infty.
\end{equation}

 By virtue of the continuity and monotonicity of $\alpha(s)$, together with the fact that $\alpha(s)<0$ on $(0,\mathfrak{S})$, we infer that $\Phi(s)$ is continuous and monotonically increasing on $(0,\mathfrak{S})$, and hence it follows from intermediary theorem that there exists an unique fixed-point $s_0\in(0,\mathfrak{S})$, such that $\Phi(s_0)=1$, in other words, $-s_0^2=\alpha(s_0)$. Taking $\lambda=s_0$, and we get $-{\lambda}^2=\alpha(\lambda)$. This proposition follows.
\endProof

It should be noted that the variational problem \eqref{3.17} achieves its minimizer on $\mathcal{A}$ according to proposition \ref{pr3.1}. Now we are on the position to show that the minimizer of problem \eqref{3.17} is also a solution to a boundary value problem equivalent to problem \eqref{1.14}-\eqref{1.15}, where $\lambda > 0$.

\begin{prop}\label{pr3.6}
If $\psi\in\mathcal{A}$ is a minimizer of problem \eqref{3.17} and denote $-\lambda^2:=E(\psi,\lambda)$, then $\psi$ is a solution to boundary value problem \eqref{1.14}-\eqref{1.15}.
\end{prop}

\Proof
For any $t,r\in\mathbb{R}$ and $\psi_0\in H_0^1\cap H^2$, define $j(t,r):=J(\psi+t\psi_0+r\psi).$ Then, $j(t,r)$ is a smooth function of $(t,r)$ with $j(0,0)=1$. Note that
\begin{align*}
\p_tj(0,0)=2\int_0^1\bar\rho(\xi^2\psi\psi_0+\psi'\psi_0')dy, \quad\p_rj(0,0)=2\int_0^1\bar\rho(\xi^2\psi^2+|\psi'|^2)dy=2\neq 0.
\end{align*}
By implicit function theorem, there exists a smooth function $r=r(t)$ near $t=0$ satisfying $r(0)=0$ and $j(t,r(t))=1$.

Since $\psi$ is the minimizer of $E(\cdot)$ on $\mathcal{A}$, one deduces that the single-variable smooth function $e(t)=E(\psi+t\psi_0+r(t)\psi)$ reaches its minimum at $t=0$, which, by Fermat's Lemma, implies that $e'(0)=0$, that is
\begin{align}\label{3.23}
  e'(0)=&2r'(0)E(\psi)+2\lambda\left[\int_0^1\mu\psi''\psi_0''dy-k_1\psi'(1)\psi'_0(1)-k_0\psi'(0)\psi'_0(0)\right]\n\\
  &-2\int_0^1g\xi^2\bar\rho'\psi\psi_0dy+2\lambda\mu\int_0^1(2\xi^2\psi'\psi'_0+\xi^4\psi\psi_0)dy=0.
\end{align}

To get $r'(0)$, we differentiate the equation $j(t,r(t))=1$ at $t=0$ and yields
\begin{align*}
2r'(0)J(\psi)+2\int_0^1\bar\rho(\xi^2\psi\psi_0+\psi'\psi_0)dy=0,
\end{align*}
which implies that
\begin{align}\label{3.24}
 r'(0)=-\int_0^1\bar\rho(\xi^2\psi\psi_0+\psi'\psi'_0)dy.
 \end{align}

Substituting (\ref{3.24}) into (\ref{3.23}) and using $E(\psi)=-\lambda^2$, one sees that
\begin{align}\label{3.25}
  -\lambda^2\int_0^1\bar\rho\left(\xi^2\psi\psi_0+\psi'\psi'_0 \right)dy
  =\lambda\left[\int_0^1\mu\psi''\psi''_0dy-(k_1\psi'(1)\psi'_0(1)+k_0\psi'(0)\psi'_0(0))\right]\nonumber\\
  +\lambda\mu\int_0^1(2\xi^2\psi'\psi'_0+\xi^4\psi\psi_0)dy-\int_0^1 g\xi^2\bar\rho'\psi\psi_0dy.
\end{align}
Taking $\psi_0\in C_0^{\infty}(0,1)$ in \eqref{3.25} infers that the minimizer $\psi$ satisfies (\ref{1.14}) in the weak sense. Then, by standard bootstrap arguments, one can deduce that $\psi$ is smooth.

Now, it remains to show that the minimizer $\psi$ satisfies boundary conditions (\ref{1.15}). In fact, since $\psi\in\mathcal{A}$, it is clear that $\psi(0)=\psi(1)=0$. Moreover, since $\psi_0\in H_0^1\cap H^2$, in view of (\ref{1.14}), applying integration by parts to (\ref{3.25}), one obtains that
\[(\mu\psi''(1)-k_1\psi'(1))\psi'_0(1)=(\mu\psi''(0)+k_0\psi'(0))\psi'_0(0).\]
As $\psi_0$ being arbitrary, we deduce that
\begin{align}\label{3.26}
\psi''(1)=\frac{k_1}{\mu}\psi'(1), \quad\psi''(0)=-\frac{k_0}{\mu}\psi'(0).
\end{align}
which implies that $\psi$ satisfies boundary conditions (\ref{1.15}).

In conclusion, $\psi$ is the solution to \eqref{1.14} with boundary condition \eqref{1.15}.
\endProof

So far, we have proved that there exists smooth solution $\psi(y)$ with the corresponding eigenvalue $\lambda>0$  for problem (\ref{1.14})-(\ref{1.15}). To ensure the validity of \rm{Fourier} synthesis in constructing exponential growing mode solution to \eqref{1.8}-\eqref{1.7}, we still need some properties of function $\lambda(\xi)$.

\begin{prop}\label{pr3.7}
The eigenvalue function $\lambda(\xi):\mathbb{A}^g\rightarrow(0,\infty)$ defined by \eqref{3.17} is continuous and bounded.
\end{prop}
\Proof
Since the definition of $-\lambda^2=\inf\limits_{\psi\in\mathcal{A}}E(\psi,\lambda)$ is equivalent to $\lambda^2=\sup\limits_{\psi\in\mathcal{A}}(-E(\psi,\lambda))$, it follows from \eqref{3.4} that
\[\lambda^2
  \leq -E(\psi)+1
  \leq g\left\|\frac{\bar\rho'}{\bar\rho}\right\|_{L^{\infty}}+\lambda C_0\left\|\frac{1}{\bar\rho}\right\|_{L^{\infty}}+1,
\]
which indicates the boundedness of $\lambda(\xi)$ on $\mathbb{A}^g$.

We now turn to the proof of the continuity, which is similar to Proposition 2.5 in \cite{JiangJN}. For the reader's convenience, we give the details here to make it more clear. For any fixed $\xi_0\in\mathbb{A}^g$, let $\xi\in \mathbb{A}^g$ with $\kappa=|\xi|^2-|\xi_0|^2$, then $|\xi|\rightarrow |\xi_0|$ as $\kappa\rightarrow 0$. The first step is to show that
\begin{align}\label{3.27}
\lim\limits_{|\xi|\rightarrow |\xi_0|}\alpha(\xi,s)=\alpha(\xi_0,s), ~~~~\text{for~any~}~s\in(0,\mathfrak{S}).
\end{align}
In view of Proposition {\ref{pr3.1}}, for any $\xi\in\mathbb{A}^g$, there exists $\psi_{\xi}\in\mathcal{A}$, such that
\begin{equation}\label{3.28}
\begin{aligned}
  \alpha(\xi,s)
  =&\int_0^1 s\mu(|\psi''_{\xi}|^2+2\xi^2|\psi'_{\xi}|^2+\xi^4\psi^2_{\xi})-s\sum_{i=0}^{1}k_i|\psi'_{\xi}(i)|^2-\int_0^1 g\xi^2\bar\rho'\psi^2_{\xi}.
\end{aligned}
\end{equation}
Substituting $|\xi|^2=|\xi_0|^2+\kappa$ into \eqref{3.28}, we get
\begin{align}\label{3.29}
   \alpha(\xi,s)
  =&\int_0^1 s\mu(|\psi''_{\xi}|^2+2\xi_0^2|\psi'_{\xi}|^2+\xi_0^4\psi^2_{\xi})
   -s\sum_{i=0}^{1}k_i|\psi'_{\xi}(i)|^2-\int_0^1 g\xi_0^2\bar\rho'\psi^2_{\xi}\nonumber\\
   &+\kappa \int_0^1 \left[s\mu(2|\psi'_{\xi}|^2+2\xi_0^2|\psi'_{\xi}|^2
     +\kappa\psi^2_{\xi})-g\bar\rho'\psi_{\xi}^2\right]\nonumber\\
  \geq& \alpha(\xi_0,s)+\kappa f(\kappa,\psi_{\xi}),
\end{align}
where
\[
  f(\kappa,\psi_{\xi})
  :=\int_0^1 \left[s\mu(2|\psi'_{\xi}|^2+2\xi_0^2|\psi'_{\xi}|^2+\kappa\psi^2_{\xi})-g\bar\rho'\psi_{\xi}^2\right].
\]
Since $\psi\in\mathcal{A}$, there exists constant $\tilde{c}$, depending on $\bar{\rho}, g, a, b$ and s, such that $|f(\kappa,\psi_{\xi})|\leq \tilde{c}$.

Similarly, we also have
\begin{align}\label{3.30}
   \alpha(\xi_0,s)\geq \alpha(\xi,s)-\kappa f(-\kappa,\psi_{\xi_0}),
\end{align}
and there exists a constant $\bar{c}$ depending on $\bar{\rho}, g, a, b$ and s, such that
\[|f(-\kappa,\psi_{\xi_0})|\leq \bar{c}.\]
Thus
\begin{equation}\label{3.31}
  \kappa f(\kappa ,\psi_{\xi})\leq\alpha(\xi,s)-\alpha(\xi_0,s)\leq\kappa f(-\kappa ,\psi_{\xi_0}).
\end{equation}
The boundedness of $f(\kappa ,\psi_{\xi})$ and $f(-\kappa ,\psi_{\xi_0})$ implies that, as $\kappa \rightarrow0$, one has
\begin{equation}\label{3.32}
  \lim\limits_{|\xi|\rightarrow|\xi_0|}\alpha(\xi,s)=\alpha(\xi_0,s),~~~~\text{for~any~}~s\in(0,\mathfrak{S}).
\end{equation}
Now, denote $\lambda(\xi,s)=\sqrt{-\alpha(\xi,s)}$, then it follows from \eqref{3.32} that
\begin{equation}\label{3.33}
  \lim\limits_{|\xi|\rightarrow|\xi_0|}\lambda(\xi,s)=\lambda(\xi_0,s)~~~~\text{for~any~}~s\in(0,\mathfrak{S}).
\end{equation}

In the second step, by virtue of proposition \ref{pr3.5}, for any $\xi\in\mathbb{A}^g$, there exists a unique $s_{\xi}\in (0,\mathfrak{S})$, such that
\begin{align}\label{3.330}
s_{\xi}=\lambda(\xi,s_{\xi})\equiv\lambda(\xi).
\end{align}
Then, the main aim of this proposition is to prove
\begin{align}\label{3.331}
\lim\limits_{|\xi|\rightarrow|\xi_0|}\lambda(\xi,s_{\xi})=\lambda(\xi_0,s_{\xi_0}).
\end{align}

In view of (\ref{3.33}), for any $\e>0$, there exists a constant $\delta>0$, as $\left|~{|\xi|-|\xi_0|}~\right|<\delta$, so that
\begin{equation}\label{3.332}
  |\lambda(\xi,s_{\xi_0})-\lambda(\xi_0,s_{\xi_0})|< \e .
\end{equation}
Moreover, since that for any $\xi\in\mathbb{A}^g$, $\alpha(\xi,s)$ is monotonically increasing on $s\in (0,\mathfrak{S})$, refer to \eqref{3.18}, one also has that $\lambda(\xi,s)$ is decreasing with respect to $s\in (0,\mathfrak{S})$.

Now, using \eqref{3.330},\eqref{3.332} and the monotonicity of $\lambda$ with respect to $s$, we are able to prove \eqref{3.331}. If $s_{\xi}\leq s_{\xi_0}$, then
\[\lambda(\xi_0,s_{\xi_0})-\e< \lambda(\xi,s_{\xi_0})\leq\lambda(\xi,s_{\xi})=s_{\xi}\leq s_{\xi_0}=\lambda(\xi_0,s_{\xi_0})<\lambda(\xi_0,s_{\xi_0})+\e.\]
On the contrary, if $s_{\xi}\geq s_{\xi_0}$, then
\[\lambda(\xi_0,s_{\xi_0})-\e<\lambda(\xi_0,s_{\xi_0})=s_{\xi_0}\leq s_{\xi}=\lambda(\xi,s_{\xi})\leq \lambda(\xi,s_{\xi_0})<\lambda(\xi_0,s_{\xi_0})+\e.\]
In conclusion, \eqref{3.331} follows from the above two inequality and the proof of this proposition is completed.
\endProof


Now, we are able to construct solutions $(\phi,\psi,\pi)$ to ODE system (\ref{1.12})-\eqref{1.13}.
\begin{prop}\label{pr3.8}
For any $\xi\in\mathbb{A}^g$, there exist solutions $(\phi,\psi,\pi)(\xi,y)$ to system (\ref{1.12})-\eqref{1.13} with the corresponding eigenvalue $\lambda(\xi)>0$. Furthermore, $(\phi,\psi,\pi)(y)\in H^k$ for any $k\in \mathbb{N}$.
\end{prop}
\Proof
In view of Proposition \ref{pr3.1} and Proposition \ref{3.6}, we have constructed solution $\psi(\xi,y)$ to problem (\ref{1.14})-(\ref{1.15}) with $\psi\in\mathcal{A}\cap H^k$, for any $k\in \mathbb{N}$. According to $(\ref{1.12})_3$ and $(\ref{1.12})_1$, we have
\begin{align*}
  \phi(\xi,y)=-\xi^{-1}\psi',~
  \pi(\xi,y)=\xi^{-1}(\lambda\bar\rho\phi+\mu\xi^2\phi-\mu\phi'')=-\xi^{-2}(\lambda\bar\rho\psi'+\mu\xi^2\psi'-\mu\psi'''),
\end{align*}
which gives rise to solutions $(\phi,\psi,\pi)$ to (\ref{1.12}), and also indicates that $(\phi,\psi,\pi)(y)\in H^k$. In addition, boundary conditions \eqref{1.13} follows directly from \eqref{1.15} and $(\ref{1.12})_3$.
\endProof

\begin{rem}\label{rem1}
In view of the definition of functional $E(\psi)$ and the definition of $\lambda(\xi)$ in \eqref{3.17}, it is clear that the function $\lambda(\xi)$ with respect to $\xi$ is an even function on $\mathbb{A}^g$. In the meantime, the associated $\psi(\xi)$ constructed in Proposition \ref{3.6} is also an even function with respect to $\xi$ defined on $\mathbb{A}^g$. Subsequently, in view of proposition \ref{pr3.8}, the corresponding function $\pi(\xi,y)$ is also an even function on $\xi\in\mathbb{A}^g$, while $\phi(\xi,y)$ is an odd function on $\xi\in\mathbb{A}^g$.
\end{rem}
The next proposition provides $H^k-$estimates for the solution $(\phi,\psi,\pi)(y)$ with $\xi$ varying. To emphasize the dependence on $\xi$, we denote the solution by $(\phi,\psi,\pi)(\xi):=(\phi,\psi,\pi)(\xi,y)$.
\begin{prop}\label{pr3.9}
For any $\xi\in\mathbb{A}^g$, let $(\phi,\psi,\pi)(\xi)$ with the corresponding function $\lambda(\xi)$ be solution to problem \eqref{1.12}-\eqref{1.13}, constructed as above. Then, for any $k\in \mathbb{N}$, there exist positive constants $A_k,B_k,C_k$ depending on $a,b,\bar\rho,\mu ,k_0,k_1,g$, such that
\begin{align*}
  \left\|\psi(\xi)\right\|_{H^k}\leq A_k,\quad
  \left\|\phi(\xi)\right\|_{H^k}\leq B_k,\quad
  \left\|\pi(\xi)\right\|_{H^k}\leq C_k.
\end{align*}
\end{prop}
\Proof
Throughout the proof, $\tilde C$ is a generic positive number depending on $a,b,\bar\rho,\mu ,k_0,k_1,g$.

Firstly, since $\psi\in\mathcal{A}$, one gets $\|\psi\|^2_{L^2}>0$. Furthermore, there exists constant $A_1$, such that
\[\|\psi(\xi)\|_{H^1}\leq A_1.\]

Secondly, Proposition \ref{pr3.5} implies that there exists $\psi(\xi)\in\mathcal{A}$, such that
\begin{align}\label{3.34}
  -\lambda^2(\xi)=\alpha(\xi,\lambda)=E(\lambda({\xi}),\psi).
\end{align}
Rewriting (\ref{3.34}), similar to \eqref{3.4}, one can deduce that
\begin{align*}
  &\lambda\mu\int_0^1|\psi''|^2dy\nonumber\\
  =&-\lambda^2+\lambda\left[\sum_{i=0}^{1}k_i|\psi'(i)|^2-\mu\int_0^1\xi^2(2|\psi'|^2+\xi^2\psi^2)dy\right]+\int_0^1g\xi^2\bar\rho'\psi^2dy\nonumber\\
  \leq&-\lambda^2+\frac{\lambda\mu}{2}\int_0^1|\psi''(\xi)|^2dy+\lambda C_0\left\|\frac{1}{\bar\rho}\right\|_{L^{\infty}(0,1)} +g\left\|\frac{\bar\rho'}{\bar\rho}\right\|_{L^{\infty}(0,1)}.
\end{align*}
According to Proposition {\ref{pr3.7}}, $\lambda$ has positive bounds from upper and lower, and hence
\begin{align}\label{3.35}
  \|\psi(\xi)\|_{H^2}\leq A_2.
\end{align}

Thirdly, equation (\ref{1.14}) can be rewritten as
\begin{align}\label{3.36}
  \psi^{(4)}(\xi)
  &=\{\lambda\mu(2\xi^2\psi''-\xi^4\psi)+g\xi^2\bar\rho'\psi-\lambda^2[\xi^2\bar\rho\psi-(\bar\rho\psi')']\}/\lambda\mu\\\nonumber
  &=\{(2\lambda\mu\xi^2+\lambda^2\bar\rho)\psi''+\lambda^2\bar\rho'\psi'+\xi^2(g\bar\rho'-\lambda\mu\xi^2-\lambda^2\bar\rho)\psi\}/\lambda\mu,
\end{align}
which infers that
 \[\|\psi^{(4)}(\xi)\|_{L^2}\leq \tilde C.\]
Applying Gagliardo-Nirenberg interpolation inequality, one has
\begin{align*}
  \left\|\psi'''\right\|_{L^2}\leq \left\|\psi''\right\|_{L^2}^{\frac{1}{2}}\left\|\psi^{(4)}\right\|_{L^2}^{\frac{1}{2}}\leq\tilde C.
\end{align*}
In conclusion, one deduce
\begin{equation}\label{3.37}
  \|\psi(\xi)\|_{H^4}\leq A_4.
\end{equation}

Differentiating equation (\ref{3.36}) with respect to $y$ and using (\ref{3.37}), we find, by induction on $k$, that $\|\psi(y)\|_{H^k}\leq A_k$, for any $k\in\mathbb{N}$. Note that $\phi,\pi$ can be expressed by $\psi$ and some constants, the rest inequalities then follows.
\endProof

In the next proposition, we will construct the exponentially growing solutions of linearized problem (\ref{1.8})-(\ref{1.7}).
\begin{prop}\label{pr3.10}
Under the assumptions of Theorem \ref{th1}, let
 \begin{align}\label{3.3701}
 \Lambda=\sup_{|\xi|\in\mathbb{A}^g\cap L^{-1}\mathbb{Z}}\lambda(\xi),
 \end{align}
 then there exist a positive constant $\Lambda^*\in(2\Lambda/3,\Lambda]$ and a real-valued solution $(\varrho,\mathbf{u},q)$ to the linearized problem (\ref{1.7})-(\ref{1.8}) defined on the horizontally periodic domain $\Omega$, such that\\
 (1) For any $k\in\mathbb{N}$,
 \begin{align}\label{3.3702}
 \|(\varrho,\mathbf{u},q)(0)\|_{H^k}<+\infty;
 \end{align}
 (2) For any $t>0$, $(\varrho,\mathbf{u},q)\in H^k$, and
 \begin{align}
 &\|(\varrho,q)(t)\|_{H^k}=e^{t\Lambda^*}\|(\varrho,q)(0)\|_{H^k},\label{3.3703}\\
 &\|u_i(t)\|_{H^k}=e^{t\Lambda^*}\|u_i(0)\|_{H^k},\ i=1,2.\label{3.3704}
 \end{align}
 (3) Moreover,
 \begin{align}\label{3.3705}
 {\rm{div}}\mathbf{u}(0)=0,\ \ \  \|u_1(0)\|_{L^2}\|u_2(0)\|_{L^2}>0
 \end{align}
 \end{prop}
\Proof
 Denote
\begin{align}\label{3.38}
  \mathbf{w}(\xi,y)=-i\phi(\xi,y)\mathbf{e}_1+\psi(\xi,y)\mathbf{e}_2,
\end{align}
where $(\phi,\psi)$ with an associated growth rate $\lambda(\xi)$ is constructed in Proposition \ref{pr3.8} for any given $\xi\in\mathbb{A}^g$.
Recalling the definition of $\Lambda$, there exist $\xi\in\mathbb{A}^g\cap L^{-1}\mathbb{Z}$, such that
\[\Lambda^*:=\lambda(\xi)=\lambda(-\xi)\in(2\Lambda/3,\Lambda].\]

In view of Remark \ref{rem1}, the following real-value functions
\begin{eqnarray}
  &&\label{3.39}
  \varrho(x,y;t)=-e^{\Lambda^*t}\bar\rho'(y)w_2(\xi,y)\left[e^{ix\xi}+e^{-ix\xi}\right],\\
  &&\label{3.40}
  \mathbf{u}(x,y;t)=\Lambda^*e^{\Lambda^*t}\left[\mathbf{w}(\xi,y)e^{ix\xi}+\mathbf{w}(-\xi,y)e^{-ix\xi}\right],\\
  &&\label{3.41}
  q(x,y;t)=e^{\Lambda^*t}\pi(\xi,y)\left[e^{ix\xi}+e^{-ix\xi}\right].
\end{eqnarray}
constitute a horizontally periodic, real-value solution to linearized problem  (\ref{1.7})-(\ref{1.8}) with
\begin{align}\label{3.42}
  \left\|(\varrho,\mathbf{u},q)(0)\right\|_{H^k}\leq M_k,~\text{and} ~\|u_2(0)\|_{L^2}>0.
\end{align}
where $M_k$ are positive numbers depending on $A_k,B_k,C_k$ in Proposition \ref{pr3.9}. In addition, $\|u_1(0)\|_{L^2}\|u_2(0)\|_{L^2}>0$ is clear, since that $\psi\in\mathcal{A}$ and $\xi\phi+\psi'=0$, and thus (\ref{3.42}) follows.

It remains to prove \eqref{3.3703} and \eqref{3.3704}. We take $q(x,y;t)$ for an example, and the others can be done similarly. Recalling the expression of $q(x,y;t)$, one sees that
\begin{align}\label{3.53}
 \|q(t)\|_{H^k}= 2\Lambda^*e^{\Lambda^* t}\|\pi(\xi,y)\cos(x\cdot\xi)\|_{H^k}=e^{\Lambda^* t}\|q(0)\|_{H^k}.
\end{align}
The proof of this proposition is completed.
\endProof

\section{Energy estimates for the nonlinear perturbed problem }
As a preparation to prove nonlinear instability in next section, referring to \cite{Guo.Tice10,JiangJN,JiangJW}, we establish some energy estimates for nonlinear perturbation equations in this section.

Suppose that $(\varrho,\mathbf{u},q)$ is a strong solution to nonlinear problem (\ref{1.6})-(\ref{1.7}). Then, for convenience, we denote
\begin{align}
  &\mathcal{E}(t):=\mathcal{E}((\varrho,\mathbf{u}))(t)=\sqrt{\|\varrho(t)\|_{H^1}^2+\|\mathbf{u}(t)\|_{H^2}^2},\label{5.1}\\
&\mathcal{E}_0:=\mathcal{E}((\varrho_0,\mathbf{u}_0))=\sqrt{\|\varrho_0\|_{H^1}^2+\|\mathbf{u}_0\|_{H^2}^2},\label{5.2}
\end{align}
and assume that there exists a constant $\delta_0\in(0,1)$ depending only on $\mu,\bar\rho,g,k_0,k_1$, such that
\begin{equation}\label{5.3}
  \mathcal{E}(t)\leq\delta_0.
\end{equation}
In this section, $C$ is denoted as a generic positive number depending only $\mu,g,\bar\rho,k_0,k_1$.

For $\mathbf{x}\in\Omega$, we define $X\in\Omega_t$ as below
\begin{align*}
  \begin{cases}
  \frac{d X(\mathbf{x},t)}{dt}=\mathbf{v}(X(\mathbf{x},t),t),\\
X(\mathbf{x},0)=\mathbf{x},
  \end{cases}
\end{align*}
where $\mathbf{v}$ is the velocity of the fluid given in \eqref{1.1}. Then, recalling \eqref{1.1}$_1$, one sees that
\[\frac{d}{dt}\rho(X(\mathbf{x},t),t)=\rho_t+\frac{dX}{dt}\cdot\nabla\rho=\rho_t+\mathbf{v}\cdot\nabla\rho=0,\]
which implies
\begin{align}\label{5.4}
  \rho(X(\mathbf{x},t),t)=\rho(X(\mathbf{x},0),0)=\rho_0(\mathbf{x}).
\end{align}
Denote that
\begin{equation}\label{5.5}
  \alpha:=\inf\limits_{\mathbf{x}\in\Omega}\{\rho_0(\mathbf{x})\}>0,~~~~~ \beta:=\sup\limits_{\mathbf{x}\in\Omega}\{\rho_0(\mathbf{x})\}<\infty.
\end{equation}
Then, for any $t\in(0,T]$, it follows from \eqref{5.4} that
\begin{equation}\label{5.6}
  0<\alpha\leq\rho(\mathbf{x},t)\leq\beta<\infty,
\end{equation}
in other words
\begin{equation}\label{5.7}
  \alpha-\bar\rho\leq\varrho(\mathbf{x},t)\leq\beta-\bar\rho,
\end{equation}
which will be used to show the boundedness of density $\rho$ and $\varrho$ in the rest of this section.

\subsection{Estimates for $\|\varrho\|_{L^2}$ and $\|\mathbf{u}\|_{L^2}$}

Multiplying \eqref{1.6}$_1$ and\eqref{1.6}$_2$ by $\varrho$ and $\mathbf{u}$ respectively,  integrating by parts over $\Omega$, and using \eqref{1.6}$_3$, we have
\begin{align}
&\frac{1}{2}\frac{d}{dt}\int\varrho^2 +\int\bar\rho'\varrho u_2 =0,\label{5.8}\\
 &\frac{1}{2}\frac{d}{dt}\int\rho|\mathbf{u}|^2 +\mu\int|\nabla \mathbf{u}|^2 -
  \int_{2\pi L\mathbb{T}}\sum_{i=0}^1k_i|u_1(x,i)|^2~dx= -\int g\varrho u_2 .\label{5.9}
\end{align}
Adding them up gives
\begin{align}\label{5.10}
  &~~~\frac{1}{2}\frac{d}{dt}\int(\varrho^2+\rho|\mathbf{u}|^2) +\mu\int|\nabla \mathbf{u}|^2 \nonumber\\
  =&-\int(\bar\rho'+g)\varrho u_2 +\int_{2\pi L\mathbb{T}}(k_1|u_1(x,1)|^2+k_0|u_1(x,0)|^2)dx:=I_1+I_2,
\end{align}

As for $I_1$, using H\"older inequality and Cauchy inequality, together with $\bar\rho(y)\in C_0^\infty(0,1)$, one gets
\begin{equation}\label{5.11}
  I_1\leq (\|\bar\rho'\|_{L^{\infty}}+g)\|\varrho\|_{L^2}\|\mathbf{u}\|_{L^2}\leq \frac{(\|\bar\rho'\|_{L^{\infty}}+g)}{2}\|(\varrho,\mathbf{u})(t)\|_{L^2}^2,
\end{equation}
To estimate $I_2$, similar to \eqref{3.4}, one also has
\begin{equation}\label{5.12}
  I_2<2C_0\|\mathbf{u}\|_{L^2}^2+\frac{\mu}{2}\|\nabla \mathbf{u}\|_{L^2}^2.
\end{equation}

Substituting (\ref{5.11}) and (\ref{5.12}) into (\ref{5.10}) and using Poincar\'e inequality \eqref{2.1}, we deduce
\begin{equation}\label{5.13}
  \frac{d}{dt}\|(\varrho,\sqrt\rho\mathbf{u})(t)\|_{L^2}^2+\mu\|\mathbf{u}\|_{H^1}^2\leq C_3\|(\varrho,\sqrt\rho\mathbf{u})(t)\|_{L^2}^2,
\end{equation}
which, together with Gronwall inequality and (\ref{5.2}), gives
\begin{equation}\label{5.14}
  \|(\varrho,\sqrt\rho\mathbf{u})(t)\|_{L^2}^2+\mu\int_0^t\|\mathbf{u}(s)\|_{H^1}^2ds\leq \beta\delta_0^2 e^{C_3t}.
\end{equation}

Furthermore, using (\ref{5.6}), one yields
\begin{equation}\label{5.15}
  \|\varrho(t)\|_{L^2}^2+\|\mathbf{u}(t)\|_{L^2}^2+\int_0^t\|\mathbf{u}(s)\|_{H^1}^2ds\leq C\delta_0^2 e^{C_3t}.
\end{equation}

\subsection{Estimates for $\|\mathbf{u}_t\|_{H^1}$ and $\|\mathbf{u}\|_{H^2}$}

First, multiplying $(\ref{1.6})_2$ by $\mathbf{u}_t$ and integrating by parts over $\Omega$, we have
\begin{align}\label{5.16}
  &~~~~\frac{1}{2}\frac{d}{dt}\int\mu|\nabla\mathbf{u}|^2~ +\int\rho |\mathbf{u}_t|^2~ \nonumber\\
  =&\int_{2\pi L\mathbb{T}}\sum_{i=0}^1k_iu_1(x,i)\p_t u_1(x,i)~dx-\int\varrho g\p_t u_2~ -\int\rho\mathbf{u}\cdot\nabla\mathbf{u}\cdot\mathbf{u}_t~ \nonumber\\
  :=&J_1+J_2+J_3.
\end{align}
Similar to \eqref{5.11}-\eqref{5.12}, we also have
\begin{align}
J_1&\leq\e\|\sqrt\rho\mathbf{u}_t\|_{L^2}^2+C_\e\|\mathbf{u}\|_{H^1}^2+\e\|\nabla\mathbf{u}_t\|_{L^2}^2,
\label{5.17}\\
J_2&\leq\int\left|\frac{\varrho g}{\sqrt\rho}\right|~|\sqrt\rho\mathbf{u}_t|~ \leq\e\|\sqrt\rho \mathbf{u}_t\|_{L^2}^2
+C_\e\|\varrho\|_{L^2}^2.\label{5.18}
\end{align}
For $J_3$, using Cauchy inequality and \eqref{5.6}, we obtain
\begin{align}\label{5.19}
  J_3\leq\|\sqrt{\rho}\mathbf{u}_t\|_{L^2}\|\sqrt\rho\mathbf{u}\cdot\nabla\mathbf{u}\|_{L^2}
  \leq\e\|\sqrt{\rho}\mathbf{u}_t\|_{L^2}^2+C_\e\|\mathbf{u}\cdot\nabla\mathbf{u}\|_{L^2}^2.
\end{align}

Second, note that $(\mathbf{u},q)$ satisfies
\begin{align}\label{5.20}
\begin{cases}
  -\mu\Delta\mathbf{u}+\nabla q=-\rho \mathbf{u}_t-\rho\mathbf{u}\cdot\nabla\mathbf{u}-\varrho g e_2, ~~&\text{in~} \Omega,\\
  {\rm div} \mathbf{u}=0, ~~&\text{in~} \Omega,\\
  u_2(x,1)=u_2(x,0)=0, ~~&\text{in~} x\in 2\pi L\mathbb{T},\\
  \p_y u_1(x,1)=\frac{k_1}{\mu}u_1(x,1), ~~&\text{in~} x\in 2\pi L\mathbb{T},\\
  \p_y u_1(x,0)=-\frac{k_0}{\mu}u_1(x,0), ~~&\text{in~} x\in 2\pi L\mathbb{T}.
\end{cases}
\end{align}
Then, it follows from the Stokes estimates, reference to \rm{Theorem A.1} in \cite{S.Ding}, that
\begin{align}\label{5.21}
\|\mathbf{u}\|_{H^2}
\lesssim \|\sqrt{\rho}\mathbf{u}_t\|_{L^2}+\|\sqrt{\rho}\mathbf{u}\cdot\nabla\mathbf{u}\|_{L^2}+g\|\varrho\|_{L^2}+\|\mathbf{u}\|_{L^2}.
\end{align}
Substituting (\ref{5.17})-(\ref{5.19}) into (\ref{5.16}) and adding the result to (\ref{5.21}), we deduce
\begin{align}\label{5.22}
  &\frac{\mu}{2}\frac{d}{dt}\|\nabla\mathbf{u}(t)\|_{L^2}^2+\frac{1}{2}\|\sqrt\rho\mathbf{u}_t\|_{L^2}^2+\frac{1}{8C_2}\|\mathbf{u}\|_{H^2}^2\nonumber\\
  \lesssim &\|(\varrho,\mathbf{u})\|_{L^2}^2+C_\e\|\mathbf{u}\|_{H^1}^2+\e\|\nabla\mathbf{u}_t\|_{L^2}^2
  +\|\mathbf{u}\cdot\nabla\mathbf{u}\|_{L^2}^2.
\end{align}
To control $\|\mathbf{u}\cdot\nabla\mathbf{u}\|_{L^2}^2$, by Lemma \ref{lemma2.3} and Lemma \ref{lemma2.4}, together with (\ref{5.3}), we get
\begin{align}\label{5.23}
  \|\mathbf{u}\cdot\nabla\mathbf{u}\|_{L^2}^2
  &\leq \|\nabla\mathbf{u}\|_{L^2}^2 \|\mathbf{u}\|_{L^\infty}^2
  \leq C(\|\mathbf{u}\|_{L^2}\|\nabla\mathbf{u}\|_{H^1})(\|\mathbf{u}\|_{L^2}\|\mathbf{u}\|_{H^2})\nonumber\\
  &\leq C\|\mathbf{u}\|_{L^2}^2\|\mathbf{u}\|_{H^2}^2\leq C\bar\delta^2\|\mathbf{u}\|_{L^2}^2.
\end{align}
Substituting (\ref{5.23}) into (\ref{5.22}) gives
\begin{align}\label{5.24}
  \mu\frac{d}{dt}\|\nabla\mathbf{u}(t)\|_{L^2}^2+\|\sqrt\rho\mathbf{u}_t\|_{L^2}^2+\|\mathbf{u}\|_{H^2}^2
  \leq C\|(\varrho,\mathbf{u})\|_{L^2}^2+C_{\e}\|\mathbf{u}\|_{H^1}^2+\e\|\nabla\mathbf{u}_t\|_{L^2}^2.
\end{align}

Third, since $\mathbf{u}=\mathbf{v}-0$ satisfies $(\ref{1.1})_2$, one has
\begin{align*}
\rho\mathbf{u}_t+\rho\mathbf{u}\cdot \nabla \mathbf{u}+\nabla p=\mu\Delta\mathbf{u}-\rho g e_2,
\end{align*}
differentiating which with respect to $t$ gives
\begin{equation}\label{5.25}
\rho\mathbf{u}_{tt}+\rho\mathbf{u}\cdot \nabla\mathbf{u}_t-\mu\Delta\mathbf{u}_t+\nabla p_t=
-\rho_t(\mathbf{u}_t+\mathbf{u}\cdot\nabla \mathbf{u}+ge_2)-\rho\mathbf{u}_t\cdot\nabla \mathbf{u}.
\end{equation}
Then, testing \eqref{5.25} by $\mathbf{u}_t$, integrating by parts over $\Omega$ and using $(\ref{1.6})$, one yields
\begin{align}\label{5.27}
  &\frac{1}{2}\frac{d}{dt}\int\rho|\mathbf{u}_{t}(t)|^2 +\mu\int |\nabla\mathbf{u}_t|^2
    -\int_{2\pi L\mathbb{T}}(k_1|\p_t\mathbf{u}_1(x,1)|^2+k_0|\p_t\mathbf{u}_1(x,0)|^2)dx\nonumber\\
  =&\int\nabla\cdot(\rho\mathbf{u})(\mathbf{u}_t+\mathbf{u}\cdot\nabla\mathbf{u}+ge_2)\cdot\mathbf{u}_t
    -\int\rho\mathbf{u}_t\cdot\nabla\mathbf{u}\cdot\mathbf{u}_t ,
\end{align}
which implies that
\begin{align}\label{5.28}
  &\frac{1}{2}\frac{d}{dt}\int\rho|\mathbf{u}_{t}(t)|^2 +\mu\int|\nabla\mathbf{u}_t|^2 \nonumber\\
  \leq&\int_{2\pi L\mathbb{T}}\sum_{i=0}^1k_i|\p_t\mathbf{u}_1(x,i)|^2dx+\int[2\rho|\mathbf{u}|~|\mathbf{u}_t|~|\nabla\mathbf{u}_t|
  +\rho|\mathbf{u}|~|\mathbf{u}_t|~|\nabla\mathbf{u}|^2+\rho|\mathbf{u}|^2~|\mathbf{u}_t|~|\nabla^2\mathbf{u}|\nonumber\\
  &+\rho|\mathbf{u}|^2~|\nabla\mathbf{u}|~|\nabla\mathbf{u}_t|
  +\rho|\mathbf{u}_t|^2|\nabla\mathbf{u}|+g\rho|\mathbf{u}|~|\nabla\mathbf{u}_t|]
  :=\sum_{i=1}^7K_j.
\end{align}

Now, we estimate each term $K_i$ above, in which \eqref{2.1}-\eqref{2.5}, H\"older inequality, Young inequality and Poincar\'e inequality may be used.
\begin{align*}
K_1&\leq 2C_0\|\mathbf{u}_t\|_{L^2}^2+\frac{\mu}{2}\|\nabla\mathbf{u}_t\|_{L^2}^2
     \leq C\|\sqrt\rho\mathbf{u}_t\|_{L^2}^2+\frac{\mu}{2}\|\nabla\mathbf{u}_t\|_{L^2}^2,
     \\
K_2&\lesssim\|\mathbf{u}\|_{L^4}\|\mathbf{u}_t\|_{L^4}\|\nabla\mathbf{u}_t\|_{L^2}
     \lesssim(\|\mathbf{u}\|_{L^2}^{1/2}\|\nabla\mathbf{u}\|_{L^2}^{1/2})
     (\|\mathbf{u}_t\|_{L^2}^{1/2}\|\nabla\mathbf{u}_t\|_{L^2}^{1/2})\|\nabla\mathbf{u}_t\|_{L^2}\\
   &\lesssim\|\nabla\mathbf{u}\|_{L^2}\|\mathbf{u}_t\|_{L^2}^{1/2}\|\nabla\mathbf{u}_t\|_{L^2}^{3/2}
     \leq C_\e\|\nabla\mathbf{u}\|_{L^2}^4\|\mathbf{u}_t\|_{L^2}^2+\e\|\nabla\mathbf{u}_t\|_{L^2}^2,
     \\
K_3&\lesssim\|\mathbf{u}\|_{L^\infty}\|\mathbf{u}_t\|_{L^2}\|\nabla\mathbf{u}\|_{L^4}^2
   \lesssim(\|\mathbf{u}\|_{L^2}^{1/2}\|\mathbf{u}\|_{H^2}^{1/2})~\|\mathbf{u}_t\|_{L^2}~
              (\|\mathbf{u}\|_{L^2}^{1/2}\|\nabla\mathbf{u}\|_{H^1}^{3/2})\\
              &\lesssim\|\mathbf{u}_t\|_{L^2}~\|\mathbf{u}\|_{L^2}~\|\mathbf{u}\|_{H^2}^2
              \lesssim\|\mathbf{u}\|_{H^1}^2\|\mathbf{u}\|_{H^2}^4+\|\mathbf{u}_t\|_{L^2}^2,
   \\
K_4&\lesssim\|\mathbf{u}\|_{L^\infty}^2\|\mathbf{u}_t\|_{L^2}\|\nabla^2\mathbf{u}\|_{L^2}
     \lesssim(\|\mathbf{u}\|_{L^2}\|\mathbf{u}\|_{H^2})\|\mathbf{u}_t\|_{L^2}\|\nabla^2\mathbf{u}\|_{L^2}
   \lesssim\|\mathbf{u}\|_{L^2}^2\|\mathbf{u}\|_{H^2}^4+\|\mathbf{u}_t\|_{L^2}^2,
   \\
K_5&\lesssim\|\mathbf{u}\|_{L^\infty}^2\|\nabla\mathbf{u}\|_{L^2}\|\nabla\mathbf{u}_t\|_{L^2}
     \lesssim(\|\mathbf{u}\|_{L^2}\|\mathbf{u}\|_{H^2})\|\nabla\mathbf{u}\|_{L^2}\|\nabla\mathbf{u}_t\|_{L^2}\\
   &\leq C_\e\|\mathbf{u}\|_{L^2}^2\|\mathbf{u}\|_{H^2}^4+\e\|\nabla\mathbf{u}_t\|_{L^2}^2,
   \\
K_6&\lesssim\|\mathbf{u}_t\|_{L^4}^2\|\nabla\mathbf{u}\|_{L^2}
     \lesssim\|\mathbf{u}_t\|_{L^2}\|\nabla\mathbf{u}_t\|_{L^2}\|\nabla\mathbf{u}\|_{L^2}
   \leq C_\e\|\mathbf{u}_t\|_{L^2}^2\|\nabla\mathbf{u}\|_{L^2}^2+\e\|\nabla\mathbf{u}_t\|_{L^2}^2,
   \\
K_7&\leq C_\e\|\mathbf{u}\|_{L^2}^2+\e\|\nabla\mathbf{u}_t\|_{L^2}^2.
\end{align*}

Putting these estimates into (\ref{5.28}) with $\e$ small sufficiently leads to
\begin{align}\label{5.29}
  &\quad\frac{1}{2}\frac{d}{dt}\|\sqrt\rho\mathbf{u}_{t}(t)\|_{L^2}^2+\frac{\mu}{4}\|\nabla\mathbf{u}_t\|_{L^2}^2\nonumber\\
  &\lesssim(1+\|\nabla\mathbf{u}\|_{L^2}^2+\|\nabla\mathbf{u}\|_{L^2}^4)\|\sqrt\rho\mathbf{u}_t\|_{L^2}^2
        +\|\mathbf{u}\|_{H^1}^2+\|\mathbf{u}\|_{L^2}^2,
\end{align}
where we have used $\|\mathbf{u}\|_{H^2}<\mathcal{E}(t)\leq\bar\delta$.

Now, choosing some constant $k_1,k_2$ large enough and $k_1>>k_2>0$, adding up $k_1\times$(\ref{5.13}), $k_2\times$(\ref{5.24}) and (\ref{5.29}), with sufficiently small $\e>0$, one can deduce that
\begin{align}\label{5.30}
  \frac{d}{dt}\|(\varrho,\sqrt\rho\mathbf{u},\nabla\mathbf{u},\sqrt\rho\mathbf{u}_t)\|_{L^2}^2+\|\mathbf{u}_t\|_{H^1}^2+\|\mathbf{u}\|_{H^2}^2
  \leq C\|(\varrho,\mathbf{u})\|_{L^2}^2.
\end{align}

In order to estimate $\|\sqrt\rho\mathbf{u}_{t}(0)\|_{L^2}^2$, we test $(\ref{1.6})_2$ by $\mathbf{u}_t$ and get
\begin{align*}
  \int\rho|\mathbf{u}_t|^2
  =&\int(-\varrho ge_2-\rho\mathbf{u}\cdot\nabla\mathbf{u}+\mu\Delta\mathbf{u})\cdot\mathbf{u}_t \\
  \leq& C\int(\varrho^2+|\mathbf{u}|^2|\nabla\mathbf{u}|^2+|\Delta\mathbf{u}|^2) +\e \int\rho|\mathbf{u}_t|^2 \\
  \leq&C\int(\varrho^2+|\Delta\mathbf{u}|^2) +\|\mathbf{u}\|_{L^\infty}^2\|\nabla\mathbf{u}\|_{L^2}^2
      +\e \int\rho|\mathbf{u}_t|^2 \nonumber\\
  \leq&C\int(\varrho^2+|\Delta\mathbf{u}|^2) +C\|\mathbf{u}\|^2_{H^2}\|\mathbf{u}\|^2_{H^2}
      +\e \int\rho|\mathbf{u}_t|^2 .
\end{align*}
Recalling $\|\mathbf{u}\|_{H^2}<\mathcal{E}(t)\leq\bar\delta$, then there holds
\begin{align*}
  \|\sqrt\rho\mathbf{u}_t\|_{L^2}^2\leq C(\|\varrho\|_{L^2}^2+\|\mathbf{u}\|_{H^2}^2)\leq C\mathcal{E}^2(t).
\end{align*}
So it follows from taking $t\to0$ that
\begin{align}\label{5.31}
  \limsup\limits_{t\to 0}\|\sqrt\rho\mathbf{u}_t\|_{L^2}^2\leq C\mathcal{E}_0^2.
\end{align}
Consequently, integrating (\ref{5.30}) with respect to the time variable on $(0,t)$ and we obtain
\begin{align}\label{5.32}
  \|(\varrho,\mathbf{u},\nabla\mathbf{u},\mathbf{u}_t)\|_{L^2}^2+\int_0^t(\|\mathbf{u}_{t}(s)\|_{H^1}^2+\|\mathbf{u}(s)\|_{H^2}^2)ds
  \leq C(\mathcal{E}_0^2+\int_0^t\|(\varrho,\mathbf{u})(s)\|_{L^2}^2ds).
\end{align}

Furthermore, recalling that $(\mathbf{u},q)$ solves the following Stokes equation
\begin{align}\label{5.33}
  -\mu\Delta\mathbf{u}+\nabla q=-\rho\mathbf{u}_t- \rho\mathbf{u}\cdot\nabla\mathbf{u}-\varrho ge_2,
\end{align}
Then, by virtue of Stokes estimate and $L^4$ estimate, we have
\begin{align}\label{5.34}
  \|\mathbf{u}\|_{H^2}^2+\|\nabla q\|_{L^2}^2
  \lesssim&\|\rho\mathbf{u}_t\|_{L^2}^2+\|\rho\mathbf{u}\cdot\nabla\mathbf{u}\|_{L^2}^2+\|\varrho\|_{L^2}^2\nonumber\\
  \lesssim&\|\mathbf{u}_t\|_{L^2}^2+\|\varrho\|_{L^2}^2+\|\mathbf{u}\|_{L^\infty}^2~\|\nabla\mathbf{u}\|_{L^2}^2)\nonumber\\
  \lesssim&\|\mathbf{u}_t\|_{L^2}^2+\|\varrho\|_{L^2}^2+\|\mathbf{u}\|_{L^2}\|\nabla\mathbf{u}\|^2_{L^2}\|\mathbf{u}\|_{H^2}\nonumber\\
  \leq&C_\e(\|\mathbf{u}_t\|_{L^2}^2+\|\varrho\|_{L^2}^2+\|\mathbf{u}\|^2_{L^2}\|\nabla\mathbf{u}\|_{L^2}^4)+\e\|\mathbf{u}\|_{H^2}^2,
\end{align}
which, together with $\mathcal{E}(t)\leq\bar\delta$ implies that
\begin{align}\label{5.35}
  \|\mathbf{u}\|_{H^2}^2+\|\nabla q\|_{L^2}^2\leq C\|(\varrho,\nabla\mathbf{u},\mathbf{u}_t)\|_{L^2}^2.
\end{align}

Finally, add (\ref{5.35}) to (\ref{5.32}), for any $t\in(0,T]$, we have,
\begin{align}\label{5.36}
  &\|\mathbf{u}\|_{H^2}^2+\|(\varrho,\mathbf{u}_t,\nabla q)\|_{L^2}^2+\int_0^t(\|\mathbf{u}_{t}(s)\|_{H^1}^2+\|\mathbf{u}(s)\|_{H^2}^2)ds\nonumber\\
  \leq& C(\mathcal{E}_0^2+\int_0^t\|(\varrho,\mathbf{u})(s)\|_{L^2}^2ds).
\end{align}

\subsection{Estimates for $\|\varrho\|_{H^1}$}

Similar to the preparation part in section 4, for any $\mathbf{x}\in\Omega$, we define the streamline function $\mathbf{X}=\mathbf{X}(\mathbf{x},t)$ by
\begin{align*}
  \begin{cases}
    \frac{d\mathbf{X}(\mathbf{x},t)}{dt}=\mathbf{u}(\mathbf{X}(\mathbf{x},t),t),\\
    \mathbf{X}(\mathbf{x},0)=\mathbf{x},\\
  \end{cases}
\end{align*}
so that
\begin{align}\label{5.37}
  \frac{d}{dt}\varrho(\mathbf{X}(\mathbf{x},t),t)
  &=\varrho_t(\mathbf{X}(\mathbf{x},t),t)+\frac{d\mathbf{X}(\mathbf{x},t)}{dt}\cdot\nabla\varrho(\mathbf{X}(\mathbf{x},t),t)\nonumber\\
  &=\varrho_t+\mathbf{u}(\mathbf{X}(\mathbf{x},t),t)\cdot\nabla\varrho(\mathbf{X}(\mathbf{x},t),t)\nonumber\\
  &=-u_2(\mathbf{X}(\mathbf{x},t),t)\bar\rho'(X_2(\mathbf{x},t)).
\end{align}
Integrating it over $(0,t)$,  we have
\begin{align}\label{5.38}
  \varrho(\mathbf{X}(\mathbf{x},t),t)=\varrho_0(\mathbf{X}(\mathbf{x},0))-\int_0^tu_2(\mathbf{X}(\mathbf{x},s),s)\bar\rho'(X_2(\mathbf{x},s))ds.
\end{align}
Then, there holds that
\begin{align}\label{5.39}
  \|\varrho(\mathbf{x},t)\|_{H^1}
  &\leq\|\varrho_0(\mathbf{x})\|_{H^1}+\|\bar\rho'\|_{L^\infty}\int_0^t\|u_2(\mathbf{x},s)\|_{H^1}ds\nonumber\\
  &\leq\|\varrho_0(\mathbf{x})\|_{H^1}+CT\sup\limits_{0\leq s\leq T}\|\mathbf{u}(s)\|_{H^1},
\end{align}

It remains to estimate $\|\mathbf{u}(s)\|_{H^1}$. In fact, applying Gronwall inequality to \eqref{5.36} yields
\begin{equation}\label{5.47}
  \sup\limits_{0\leq s\leq T}\|\mathbf{u}(s)\|^2_{H^1}\leq C(T)\mathcal{E}^2_0.
\end{equation}
Then, one obtains
\begin{align}\label{5.40}
  \|\varrho(t)\|_{H^1}\leq C(T)\mathcal{E}^2_0.
\end{align}

In conclusion, it follows from (\ref{5.36}) and (\ref{5.40}) that
\begin{align}\label{5.48}
  \mathcal{E}^2(t)+\|(\mathbf{u}_t,\nabla q)\|_{L^2}^2+\int_0^t(\|\mathbf{u}_{t}(s)\|_{H^1}^2+\|\mathbf{u}(s)\|_{H^2}^2)ds
  \leq C(T)\mathcal{E}_0^2.
\end{align}

Compared to classical well-posedness results of nonhomogeneous Navier-Stokes equations in \cite{J.Choe,J.Simon}, the boundary integral terms produced by the Navier-slip boundary condition in this paper bring some new difficulties in the energy estimates. However, these difficulties can be overcome by using the technique applied in \cite{Li}. Thus, based on estimate \eqref{5.48}, the global well-posedness of the nonlinear system \eqref{1.6}-\eqref{1.7} follows. We state this result below without the proof.
\begin{prop}\label{pr5.1}
Suppose that the steady state satisfies \eqref{1.3}. Then for any given initial data $(\varrho_0,\mathbf{u}_0)\in(H^1\cap L^\infty)\times H^2$ satisfying \eqref{5.2}, \eqref{5.5}, $\nabla\cdot\mathbf{u}_0=0$, and also being compatible with the boundary conditions \eqref{1.2}, the nonlinear problem \eqref{1.6}- \eqref{1.7} has a global strong solution $(\varrho,\mathbf{u},\nabla q)\in C([0,T];H^1\times H^2\times L^2)$, such that
\begin{align}\label{5.49}
 \mathcal{E}^2(t)+\|(\mathbf{u}_t,\nabla q)\|_{L^2}^2+\int_0^t(\|\mathbf{u}_{t}(s)\|_{H^1}^2+\|\mathbf{u}(s)\|_{H^2}^2)ds
  \leq C_4(T)\mathcal{E}_0^2,
\end{align}
where the positive number $C_4(T)$ depends only on $g,\mu,\alpha,\beta,\bar\rho$ and also $T$.
\end{prop}

\section{The proof of nonlinear instability}
In this section, we will prove the nonlinear instability by the bootstrap argument proposed by Y. Guo et al. in \cite{Guo.S}. To be precise, we will show that there exists a constant $\varepsilon > 0 $ such that for any $\delta_0$, though small enough, and the initial data is smaller then $\delta_0$ in some suitable sense, the nonlinear system \eqref{1.6}-\eqref{1.7} admits a strong solution ${\mathbf u}^\delta$ and an escape time $T^\delta > 0$ such that $\|{\mathbf u}^\delta(T^\delta)\|_{L^2} > \varepsilon$.

To this end, we first give the following elementary inequality, which will be used in the proof of Proposition \ref{prop6.2}.
\begin{lem}\label{lemma6.1}
Let ${\mathbf w}\in H_\sigma^1\cap H^2$, where $H_\sigma^1=\{{\mathbf w}\in H^1| {\mathrm{div}}{\mathbf w}=0, w_2=0~{\mathrm{on}}~\Sigma_1\bigcup\Sigma_0\}$, then it holds that
\begin{equation}\label{6.1}
\begin{aligned}
  \int g\bar{\rho}'|{w}_2|^2+\Lambda\sum_{i=0}^1 k_i\int_{2\pi L\mathbb{T}}|w_1(x,i)|^2{\mathrm{d}}x
  \leq \Lambda^2{\int\bar{\rho}|{\mathbf w}|^2}+\Lambda\mu\int{|\nabla {\mathbf w}|^2}.
\end{aligned}
\end{equation}
\end{lem}
\Proof
Since the proof is similar and simpler to the Step 1 of Lemma 4.1 in \cite{JJW}, we omit the details here.
\endProof

According to Proposition \ref{pr3.10},  there exists constant $\Lambda^*\in(2\Lambda/3,\Lambda]$, such that the exponentially increasing functions
\begin{equation}\label{6.10}
  \left(\varrho^l,{\mathbf u}^l\right)
  =e^{\Lambda^* t}\left(\bar\varrho_0,\bar{\mathbf u}_0\right)\in H^2 \times (H_\sigma^1\cap H^2)~~{\mbox{for each}}~~ t > 0
\end{equation}
satisfy the linearized system \eqref{1.8}-\eqref{1.7} with an associated pressure $q^l=e^{\Lambda^* t}\bar{q}_0$, where $\bar{q}_0\in H^1$, and
$(\bar{\varrho}_0,\bar{ u}_0)\in H^2\times( H_\sigma^1\cap H^2)$ satisfy
\begin{equation} \label{6.11}
  \|\bar{\varrho}_0\|_{L^2}\|{\bar u}_{02}\|_{L^2}\|{\bar{ u}}_{01}\|_{L^2}>0,\quad
  {\mathcal{E}}((\bar{\varrho}_0,\bar{\mathbf u}_0))=\sqrt{\|\bar{\varrho}_0\|_{H^1}^2+\|\bar{\mathbf u}_0\|_{H^2}^2}=1.
\end{equation}
Here $\bar{u}_{0i}$  stands for the $i$-th component of $\bar{{\mathbf u}}_0$ for $i=1,2$.

Denote
$(\varrho_0^\delta,{\mathbf u}_0^\delta):=\delta (\bar{\varrho}_0,\bar{{\mathbf u}}_0)$, and $C_5:=\|(\bar{\varrho}_0,\bar{{\mathbf u}}_0 )\|_{L^2}$.
Keeping in mind that \[\inf_{ x\in\Omega}\{\bar{\rho}( x)\}>0\] and the embedding $H^2\hookrightarrow L^\infty$, we can choose a sufficiently small $\delta_1\in (0,1)$, such that
\begin{equation*}
  \frac{\inf_{ x\in\Omega}\{\bar{\rho}( x)\}}{2}\leq
  \inf_{ x\in\Omega}\{\varrho_0^\delta( x)+\bar{\rho}( x)\} \mbox{ for any }\delta\in (0,\delta_1).
\end{equation*}
Hence, by virtue of Proposition \ref{pr5.1}, the perturbed problem \eqref{1.6}--\eqref{1.7} admits a strong solution
$(\varrho^\delta, \mathbf u^\delta)\in C^0([0,T^{\max}),H^1\times H^2)$ with an associated pressure $\nabla q^\delta\in C^0([0,{T^{\max}}),L^2)$, satisfying
the initial data $(\varrho_0^\delta, {\mathbf u}_0^\delta)$ with ${\mathcal{E}}((\varrho_0^\delta, {\mathbf u}_0^\delta))=\delta$. Moreover, we have
\begin{equation}\label{6.12}
  0<\frac{\inf_{ x\in\Omega}\{\bar{\rho}( x)\}}{2}\leq\inf_{ x\in\Omega}\{\varrho^\delta (t,x)+\bar{\rho}\}
\end{equation}
\begin{equation}\label{6.13}
  \sup_{ x\in\Omega}\{\varrho^\delta(t,x)+\bar{\rho}\}\leq \sup_{ x\in\Omega}\{\bar{\varrho}_0( x)+\bar{\rho}\}
  \leq C_6\|\bar{\varrho}_0\|_{H^2}+\|\bar{\rho}\|_{L^\infty}
\end{equation} for any $t\in [0,T^{\max})$, where $C_6$ is the constant from the embedding $H^2\hookrightarrow {L^\infty}$.

Now we choose the value of $\delta_0\in(0,1)$ as small as \eqref{5.3}. Let $\sigma=\min \{\delta_0,\delta_1,\varepsilon_0\}$, and $\delta\in (0,\sigma)$, define
\begin{equation}\label{6.14}
  T^{\delta}:=\frac{1}{\Lambda^*}{\mathrm{ln}}\frac{2\varepsilon_0}{\delta}>0,\quad\mbox{i.e.,}\;\delta e^{\Lambda^* T^\delta}=2\varepsilon_0,
\end{equation}
where $\varepsilon_0$ is a constant independent of $\delta$, satisfying $\varepsilon_0\in (0,1)$, which will be defined in \eqref{6.37}.
\begin{equation*}
  T^*:=\sup\left\{t\in I_{T^{\max}}\left|~{\mathcal{E}}((\varrho^\delta,{\mathbf u}^\delta )(t))\leq {\delta_0}\right.\right\}>0,
\end{equation*}
\begin{equation*}
  T^{**}:=\sup\left\{t\in I_{T^{\max}}\left|~\left\|\left(\varrho^\delta,{\mathbf u}^\delta\right)(t)\right\|_{{L}^2 }\leq 2\delta C_5e^{\Lambda^* t}\right\}>0\right..
\end{equation*}
Then $T^*$ and $T^{**}$ may be finite, and furthermore,
 \begin{eqnarray}
   \label{6.15}
   &&{\mathcal{E}}(\left(\varrho^\delta, {\mathbf u}^\delta\right)(T^*))={\delta_0},\quad\mbox{ if }T^*<\infty ,\\
   \label{6.16}
   &&\left\|\left(\varrho^\delta, {\mathbf u}^\delta\right)(T^{**})\right\|_{{L}^2 }=2\delta C_5e^{\Lambda T^{**}},\quad
   \mbox{ if }T^{**}<T^{\max}.
\end{eqnarray}

Now, we denote ${T}_{\min}:= \min\{T^\delta,T^*,T^{**}\}$, then for all $t\in
\bar{I}_{ {T}_{\min}}$, we deduce from the estimate \eqref{5.49} and the definitions of $T^*$ and $T^{**}$ that
\begin{align} \label{6.17}
  {\mathcal{E}}^2\big( (\varrho^\delta, {\mathbf u}^\delta)(t)\big) +\|\mathbf u_{t}^\delta(t)\|_{L^2 }^2+\int_0^t
  \|\nabla\mathbf u_\tau^\delta\|_{L^2 }^2\mathrm{d}\tau
  \leq&~ C_4 \delta^2 {\mathcal{E}}^2(\left(\bar{\varrho}_0,\bar{\mathbf u}_0 \right))
  \leq~  C_4\delta^2e^{2\Lambda^* t},
\end{align}
where $C_4$ is independent of $\delta$.

Let $(\varrho^{\mathrm{d}}, {\mathbf u}^{\mathrm{d}})=(\varrho^{\delta},{\mathbf u}^{\delta})-\delta(\varrho^l, {\mathbf u}^l)$. Noting that
$(\varrho^\mathrm{a},{\mathbf u}^{\mathrm{a}}):= \delta(\varrho^l,{\mathbf u}^l)\in C^0([0,+\infty),H^1\times H^2 )$ is also a linear solution to \eqref{1.8}--\eqref{1.7} with the initial data
$(\varrho_0^\delta, \mathbf u_0^\delta)\in H^1 \times H^2 $ and an associated pressure $q^a=\delta q^l\in
 C^0([0,+\infty),L^2 )$, we find that $(\varrho^{\mathrm{d}}, {\mathbf u}^{\mathrm{d}})$
satisfies the following error equations:
\begin{equation}\label{6.18}
\left\{\begin{array}{ll}
  \varrho_t^{\mathrm{d}}+\bar{\rho}'{u}_2^{\mathrm{d}}= -{{\mathbf u}}^{\delta}\cdot \nabla\varrho^{\delta}, \\[1mm]
  (\varrho^\delta+\bar{\rho}) \mathbf u_t^{\mathrm{d}} -\mu \Delta \mathbf u^{\mathrm{d}}+\nabla q^{\mathrm{d}}
  = {\mathbf f}^\delta-g\varrho^{\mathrm{d}}{\mathbf e}_2,\\[1mm]
 \mathrm{div}\mathbf u^{\mathrm{d}}={0},
\end{array}\right.
\end{equation}
where
$$q^{\mathrm{d}}:=q^\delta-q^\mathrm{a}\in C^0(\bar{I}_{ {T}_{\min}},H^1)\mbox{ and   }
 {\mathbf f}^\delta:=-( \varrho^{\delta}+\bar{\rho}) \mathbf u^{\delta}\cdot\nabla \mathbf u^{\delta}-\varrho^\delta\mathbf u^{\mathrm{a}}_t.$$
The initial and boundary conditions become
\begin{align*}
\begin{cases}
  (\varrho^{\mathrm{d}},{{\mathbf u}}^{\mathrm{d}})|_{t=0}= 0,\\
  u_2^{\mathrm{d}}(x,0)=u_2^{\mathrm{d}}(x,1)=0,\\
  \p_y u_1^{\mathrm{d}}(x,1)=\frac{k_1}{\mu}u_1^{\mathrm{d}}(x,1),\\
  \p_y u_1^{\mathrm{d}}(x,0)=-\frac{k_0}{\mu}u_1^{\mathrm{d}}(x,0),
\end{cases}
\end{align*}
with compatibility conditions read as
$$u_2^\mathrm{d}(x,0)|_{t=0}=u_2^\mathrm{d}(x,1)|_{t=0}=0, \quad\mathrm{div} \mathbf u^\mathrm{d}|_{t=0}=0,$$

In the following, we will establish the error estimate for $(\varrho^{\mathrm{d}},\mathbf u^{\mathrm{d}})$ in $L^2$-norm.
\begin{prop}\label{prop6.2}
There is a constant $C_7$, such that for all $t\in \bar{I}_{{T}_{\min}}$,
\begin{equation}\label{6.19}
\begin{aligned}
  \| (\varrho^{\mathrm{d}},\mathbf u^{\mathrm{d}})(t)\|^2_{L^2 } \leq C_7\delta^3e^{3\Lambda^* t}.
\end{aligned}
\end{equation}
\end{prop}
\Proof
Recalling that
$(\varrho^{\mathrm{d}}, {\mathbf u}^{\mathrm{d}})=(\varrho^\delta,\mathbf u^\delta)-(\varrho^\mathrm a,\mathbf u^\mathrm a)$, in view of the regularity of $(\varrho^\delta,\mathbf u^\delta)$ and $(\varrho^\mathrm{a},\mathbf u^{\mathrm{a}})$, we can deduce from \eqref{6.18}$_2$ that for a.e. $t\in {I}_{T_{\min}}$,
\begin{equation}\label{6.20}
\begin{aligned}
  &\frac{\mathrm{d}}{\mathrm{d}t}\int (\varrho^\delta+\bar{\rho})| \mathbf u_t^{\mathrm{d}}|^2
  =2<((\varrho^\delta+\bar{\rho}){\mathbf u}_t^{\mathrm{d}})_t,{\mathbf u}_t^{\mathrm{d}}>-\int \varrho^\delta_t|\mathbf u_t^{\mathrm{d}}|^2\\
  =&2\int ({\mathbf f}_t^\delta-g\varrho^{\mathrm{d}}_t \mathbf e_2) \mathbf u^{\mathrm{d}}_t-2(\mu\int|\nabla\mathbf u^{\mathrm{d}}_t|^2-\int_{2\pi L\mathbb{T}}k_i|\p_t u_1^d(x,i)|^2\mathrm{d}x) -\int \varrho^\delta_t|\mathbf u_t^{\mathrm{d}}|^2,
\end{aligned}\end{equation}
and $\|\sqrt{\varrho^\delta+\bar{\rho}}\mathbf u_t^{\mathrm{d}}\|_{L^2}\in C^0(\bar{I}_{T_{\min}})$, referring to \cite[Remark 6]{Y.Cho}.
Noting that
\begin{equation}\label{6.201}
\frac{\mathrm{d}}{\mathrm{d}t}\int \bar{\rho}'|{u}_2^{\mathrm{d}}|^2=2\int \bar{\rho}'{u}_2^{\mathrm{d}}\partial_t u_2^{\mathrm{d}},
\end{equation}
thus, adding up the equality \eqref{6.20} and \eqref{6.201}, using \eqref{6.18}$_1$, it gives
\begin{equation}\label{6.21}
\begin{aligned}
  &\frac{\mathrm{d}}{\mathrm{d}t}\int \left[(\varrho^\delta+\bar{\rho})|\mathbf u_t^{\mathrm{d}}|^2
  -g\bar{\rho}'|{u}_2^{\mathrm{d}}|^2\right]+2(\mu\int|\nabla\mathbf u^{\mathrm{d}}_t|^2-\int_{2\pi L\mathbb{T}}k_i|\p_t u_1^d(x,i)|^2\mathrm{d}x)\\
  =&~\int\left (2{\mathbf f}^\delta_t+2g\mathbf u^\delta\cdot\nabla \varrho^\delta \mathbf e_2- \varrho^\delta_t\mathbf u_t^{\mathrm{d}}\right)\cdot \mathbf u_t^{\mathrm{d}},
\end{aligned}\end{equation}

Integrating (\ref{6.21}) with respect to time variable from $0$ to $t$, we get
\begin{align}\label{6.22}
  &\|\sqrt{\varrho^\delta+\bar{\rho}}\mathbf u_t^\mathrm{d}(t)\|^2_{L^2}+2\int_0^t(\mu\|\nabla\mathbf u_\tau ^\mathrm{d}\|^2_{L^2 }-\int_{2\pi L\mathbb{T}}k_i|\p_\tau u_1^d(x,i)|^2\mathrm{d}x)\mathrm{d}\tau\nonumber\\
  =&~\int g\bar{\rho}'|{u}^{\mathrm{d}}_2(t)|^2 + R_1+R_2(t),
\end{align}
where
\begin{equation*}
  R_1=\left[\int(\varrho^\delta+\bar{\rho})| \mathbf u_t^{\mathrm{d}}|^2\mathrm{d} x\right]_{t=0}
\end{equation*}
and
\begin{equation*}
  R_2(t)=\int_0^t
  \left(2{\mathbf f}_\tau+2 g\mathbf u^\delta\cdot\nabla \varrho^\delta\mathbf e_2- \varrho^\delta_\tau\mathbf u_\tau^{\mathrm{d}}\right)\cdot \mathbf u_\tau^{\mathrm{d}} \mathrm{d}\tau.
\end{equation*}
The estimate of the above two terms $R_1$ and $R_2(t)$ follows from \cite{JJW}. For convenience, we state the conclusion without proofs
\begin{equation}\label{6.23}
  R_1+R_2(t)\lesssim \delta^3 e^{3\Lambda^* t},
\end{equation}
which, together with \eqref{6.22}, yields that
\begin{align}\label{6.24}
  &\|\sqrt{\varrho^\delta+\bar{\rho}}\mathbf u_t^\mathrm{d}(t)\|^2_{L^2}+2\int_0^t(\mu\|\nabla\mathbf u_\tau ^\mathrm{d}\|^2_{L^2}-\int_{2\pi L\mathbb T}k_i|\p_\tau u_1^\mathrm{d}(x,i)|^2\mathrm{d}x)\mathrm{d}\tau\nonumber\\
  \leq& \int g\bar{\rho}'|{ u}^{\mathrm{d}}_2|^2+C\delta^3 e^{3\Lambda^* t}.
 \end{align}

In addition, it follows from Lemma \ref{lemma6.1} that
\begin{align}\label{6.25}
  \int g\bar{\rho}'|{u}^{\mathrm{d}}_2|^2
  \leq~\Lambda^2{\int(\varrho^\delta+\bar{\rho})|{\mathbf u}^{\mathrm{d}}|^2}
   +\Lambda\mu\int{|\nabla {\mathbf u}^{\mathrm{d}}|^2}-\Lambda\sum_{i=0}^1 k_i\int_{2\pi L\mathbb{T}}|u_1^{\mathrm{d}}(x,i)|^2{\mathrm{d}}x.
\end{align}

Recalling that $\mathbf u^\mathrm{d}\in C^0(\bar{I}_{T_{\min}}, H^2)$ and $\nabla\mathbf u^{\mathrm{d}}(t)|_{t=0}=0$, using Newton-Leibniz's formula and Cauchy inequality, we rewrite and estimate the last two terms in the right hand side as follows:
\begin{align}\label{6.27}
  &\Lambda(\mu\int{|\nabla {\mathbf u}^{\mathrm{d}}|^2}-\sum_{i=0}^1 k_i\int_{2\pi L\mathbb{T}}|u_1^{\mathrm{d}}(x,i)|^2{\mathrm{d}}x)\nonumber\\
  \leq&~\Lambda^2 \int_0^t(\mu\|\nabla\mathbf u^{\mathrm{d}}\|^2_{L^2}-\sum_{i=0}^1 \int_{2\pi L\mathbb{T}}k_i| u_1^{\mathrm d}(x,i)|^2\mathrm{d}x)\mathrm{d}\tau\nonumber\\
  &~+\int_0^t(\mu\|\nabla\mathbf u_\tau ^{\mathrm{d}}\|^2_{L^2}-\sum_{i=0}^1 \int_{2\pi L\mathbb{T}}k_i|\p_\tau u_1^{\mathrm d}(x,i)|^2\mathrm{d}x)\mathrm{d}\tau\nonumber\\
  &~+\int_0^t[\sum_{i=0}^1 \int_{2\pi L\mathbb{T}}k_i(\Lambda u_1^{\mathrm{d}}(x,i)-\p_\tau u_1^{\mathrm d}(x,i))^2\mathrm{d}x]\mathrm{d}\tau.
\end{align}
Putting the above three inequalities together gives
\begin{align}\label{6.28}
  &\|\sqrt{\varrho^\delta+\bar{\rho}}\mathbf u_t^\mathrm{d}(t)\|^2_{L^2}+\frac{1}{2}\int_0^t(\mu\|\nabla\mathbf u_\tau ^\mathrm{d}\|^2_{L^2}-\sum_{i=0}^1 \int_{2\pi L\mathbb{T}}k_i|\p_\tau u_1^{\mathrm d}(x,i)|^2\mathrm{d}x)\mathrm{d}\tau\nonumber\\
  &+\frac{1}{2}\Lambda(\mu\int{|\nabla {\mathbf u}^{\mathrm{d}}|^2}-\sum_{i=0}^1 k_i\int_{2\pi L\mathbb{T}}|u_1^{\mathrm{d}}(x,i)|^2{\mathrm{d}}x)\nonumber\\
  \leq&~~\Lambda^2{\int(\varrho^\delta+\bar{\rho})|{\mathbf u}^{\mathrm{d}}|^2}+C\delta^3 e^{3\Lambda^* t}+\frac{3}{2}\Lambda^2 \int_0^t(\mu\|\nabla\mathbf u^{\mathrm{d}}\|^2_{L^2}-\sum_{i=0}^1 \int_{2\pi L\mathbb{T}}k_i| u_1^{\mathrm d}(x,i)|^2\mathrm{d}x)\mathrm{d}\tau\nonumber\\
  &~+\frac{3}{2}\int_0^t[\sum_{i=0}^1 \int_{2\pi L\mathbb{T}}k_i(\Lambda u_1^{\mathrm{d}}(x,i)-\p_\tau u_1^{\mathrm d}(x,i))^2\mathrm{d}x]\mathrm{d}\tau.
\end{align}

On the other hand, by virtue of Cauchy inequality, we get
\begin{equation}\begin{aligned}\label{6.29}
  \frac{\mathrm{d}}{\mathrm{d}t}\|\sqrt{\varrho^\delta+\bar{\rho}}\mathbf u^\mathrm{d} \|^2_{L^2}
  =&2\int(\varrho^\delta+\bar{\rho})\mathbf u^\mathrm{d} \cdot \mathbf u^\mathrm{d}_t+\int\varrho^\delta_t |\mathbf u^\mathrm{d}|^2\\
  \leq&\frac{1}{\Lambda}\|\sqrt{(\varrho^\delta+\bar{\rho})} \mathbf u_t^\mathrm{d} \|^2_{L^2}
  +\Lambda\|\sqrt{\varrho^\delta+ \bar{\rho}}\mathbf u^\mathrm{d} \|^2_{L^2}+\int\varrho^\delta_t |\mathbf u^\mathrm{d}|^2.
\end{aligned}\end{equation}
Utilizing \eqref{6.10},\eqref{6.13},\eqref{6.17} and embedding inequality, the last term can be estimated as
\begin{equation}\label{6.30}
\begin{aligned}
  \int\varrho^\delta_t |\mathbf u^\mathrm{d}|^2
  =-\int({\mathbf u}^{\delta}\cdot \nabla\varrho^{\delta}+\bar{\rho}'{u}_2^{\delta}) |\mathbf u^\mathrm{d}|^2
  =\int(2\varrho^{\delta}{\mathbf u}^{\delta}\cdot \nabla\mathbf u^\mathrm{d}-\bar{\rho}'{u}_2^{\delta}\mathbf u^\mathrm{d})\cdot \mathbf u^\mathrm{d}
  \lesssim \delta^3 e^{3\Lambda^* t}.
\end{aligned}
\end{equation}
Furthermore, integrating \eqref{6.29} in time from 0 to t, together with \eqref{6.28} and using Young inequality, we have
\begin{align}\label{6.31}
  &\|\sqrt{\varrho^\delta+\bar{\rho}}\mathbf u_t^\mathrm{d}\|^2_{L^2}+\Lambda\mu\|\nabla {\mathbf{u}}^{\mathrm{d}}\|^2_{L^2}+\int_0^t\mu\|\nabla\mathbf u_\tau ^\mathrm{d}\|^2_{L^2}\mathrm{d}\tau\nonumber\\
  \leq&~\Lambda\int_0^t(\|\sqrt{\varrho^\delta+\bar{\rho}}\mathbf u_\tau^\mathrm{d}\|^2_{L^2}+\Lambda\mu\|\nabla {\mathbf{u}}^{\mathrm{d}}\|^2_{L^2})\mathrm{d}\tau+\Lambda\|\sqrt{\varrho^\delta+\bar{\rho}}\mathbf u^\mathrm{d} \|^2_{L^2}+C\delta^3 e^{3\Lambda^* t}.
\end{align}

Summing up the previous three estimates together yields

\begin{align}\label{6.32}
  &\frac{\mathrm{d}}{\mathrm{d}t} \|\sqrt{\varrho^\delta+\bar{\rho}}\mathbf u^\mathrm{d}(t)\|^2_{L^2}+(\|\sqrt{\varrho^\delta+\bar{\rho}}\mathbf u_t^\mathrm{d}(t)\|^2_{L^2}+ \Lambda\mu\|\nabla \mathbf u^\mathrm{d}(t)\|_{L^2}^2)\nonumber\\
  \leq&~\Lambda\left[\|\sqrt{\varrho^\delta+\bar{\rho}}\mathbf u^\mathrm{d}(t)\|^2_{L^2}+\int_0^t(\|\sqrt{\varrho^\delta+\bar{\rho}}\mathbf u_\tau^\mathrm{d}\|^2_{L^2}+ \Lambda\mu\|\nabla \mathbf u^\mathrm{d}\|_{L^2}^2)\mathrm{d}\tau\right]+C\delta^3 e^{3\Lambda^* t}.
\end{align}

Therefore, applying Gronwall's inequality to \eqref{6.32}, one obtains
\begin{equation}\label{6.33}
\begin{aligned}
  \|\sqrt{\varrho^\delta+\bar{\rho}}\mathbf u^\mathrm{d}(t)\|^2_{L^2 }+\int_0^t(\|\sqrt{\varrho^\delta+\bar{\rho}}\mathbf u_\tau^\mathrm{d}\|^2_{L^2 }+ \|\nabla \mathbf u^\mathrm{d}\|_{L^2 }^2)\mathrm{d}\tau\lesssim C\delta^3 e^{3\Lambda^* t},
\end{aligned}
\end{equation}
for all $t\leq \bar{I}_{T_{\min}}$, where the constant $C$ depends on $k_0,k_1,\mu,\bar\rho,T_{min},\Lambda$. Together with \eqref{6.13} and \eqref{6.31}, we deduce that
\begin{eqnarray}\label{6.34}
 \| \mathbf u^{\mathrm{d}}(t)\|_{H^1 }^2+\int_0^t(\| \mathbf u_\tau^{\mathrm{d}}\|^2_{L^2 }+
 \|\nabla\mathbf u^{\mathrm{d}}\|^2_{L^2 })\mathrm{d}\tau \lesssim \delta^3e^{3\Lambda^* t}.
\end{eqnarray}

Finally, using the estimates \eqref{6.17}, \eqref{6.34} and embedding inequality, we can deduce from the equation \eqref{6.18}$_1$ that
\begin{equation}\begin{aligned}\label{6.36}
  \|\varrho^{\mathrm{d}}(t)\|_{L^2 }
  \leq &\int_0^t\|\varrho^{\mathrm{d}}_\tau\|_{L^2 }\mathrm{d}\tau
  \lesssim \int_0^t(\| \mathbf u^\mathrm{d} \|_{H^1 }+\|{{ \mathbf u}}^{\delta}\cdot \nabla\varrho^{\delta}\|_{L^2 })\mathrm{d}\tau \\
  \lesssim & \int_0^t(\delta^\frac{3}{2}e^{\frac{3\Lambda^*}{2}\tau}+
  \delta^2 e^{2\Lambda^*\tau})\mathrm{d}\tau\lesssim \delta^\frac{3}{2}e^{\frac{3\Lambda^*}{2} t},
\end{aligned}\end{equation}
in which we have used the following fact
\begin{equation}\label{6.35}
  \delta e^{\Lambda^* t}\leq \delta e^{\Lambda^* T^\delta}=\e_0<1~~~\textrm {for~any }t\in \bar{I}_{T_{\min}}.
\end{equation}
Therefore, consolidating \eqref{6.36} with \eqref{6.33}, the desired estimate \eqref{6.19} follows. This completes the proof of this proposition.
\endProof

\noindent{\bf Proof of Theorem 1.2.~}Now, with \eqref{6.10},\eqref{6.14},\eqref{6.17} and \eqref{6.19} in hands, referring to \cite{JJW} for details, we can conclude that
\begin{equation}\label{6.37}
  T^\delta=T_{\min},~~ {\mathrm{provided}} ~~\varepsilon_0=\min\left\{\frac{{\delta_0}}{4\sqrt C_4},\frac{C_5^2}{2C_7^2},\frac{m_0^2}{8C_7} \right\},
\end{equation}
 where we have defined that
$m_0=: \min\{ \|\bar{\varrho}_0\|_{L^2},\|\bar{u}_{02}\|_{L^2},\|\bar{u}_{01}\|_{L^2}\}>0$ due to \eqref{6.11}.

Since $T^\delta=T_{\min}$, \eqref{6.19} holds for $t=T^\delta$. Therefore,
we can use \eqref{6.19} and \eqref{6.37} with $t=T^\delta$ to obtain that
\begin{equation*}
 \begin{aligned}
  \|\varrho^{\delta}(T^\delta)\|_{L^2 }\geq &
  \|\varrho^{\mathrm{a}}_{\delta}(T^{\delta})\|_{L^2 }-\|\varrho^{\mathrm{d}}(T^{\delta})\|_{L^2 }
  \geq  \delta e^{\Lambda^* T^\delta}\|\bar{\varrho}_{0}\|_{L^2 } -\sqrt{C_7}\delta^{3/2}e^{3\Lambda^* T^{\delta}/2}\\
  \geq & 2\varepsilon_0\|\bar{\varrho}_{0}\|_{L^2 } -2^{3/2}\sqrt{C_7}\varepsilon_0^{3/2}
  \geq 2m_0\varepsilon_0 -2^{3/2}\sqrt{C_7}\varepsilon_0^{3/2} \geq m_0\varepsilon_0.
 \end{aligned}
\end{equation*}
 Similarly, it is easy to get that
\begin{equation*}\begin{aligned}
  \|u_i^{\delta}(T^\delta)\|_{L^2 }
  \geq  2m_0\varepsilon_0 -2^{3/2}\sqrt{C_7}\varepsilon_0^{3/2}\geq m_0\varepsilon_0,
\end{aligned}\end{equation*}
 where $u^{\delta}_{i}(T^{\delta})$ denotes the $i$-th component of
$ u^{\delta}(T^{\delta})$ for $i=1$, $2$.
This completes the proof of Theorem \ref{th2} by defining $\varepsilon :=m_0\varepsilon_0$.
\endProof

\acknowledgment
Ding's research is supported by the National Natural Science Foundation of China (No.11371152, No.11571117, No.11771155 and No.11871005) and Natural Science Foundation of Guangdong Province, China(No.2017A030313003). Li's research is supported by the National Natural Science Foundation of China (No.11901399) and the Natural Science Foundation of Shenzhen University (No.2019084).

This work does not have any conflicts of interest.


{\it dingsj@scnu.edu.cn(S. Ding), zhijunji@m.scnu.edu.cn(Z. Ji), quanrong\_li@szu.edu.cn(Q. Li)}
\end{document}